\documentclass{birkjour}
\usepackage{amsmath}
\usepackage{amsfonts}
\usepackage{amssymb}
\usepackage{dsfont}
\usepackage{newlfont}
\usepackage{mathtools}
\usepackage{stmaryrd}
\usepackage{enumitem}
\usepackage{graphicx}
\usepackage{array}
\usepackage{multirow}
\usepackage{ulem}
\usepackage[page,toc]{appendix}
 \newtheorem{theo}{Theorem}[section]
 \newtheorem{propos}[theo]{Proposition}
 \newtheorem{lemma}[theo]{Lemma}
 
 \newtheorem{coro}[theo]{Corollary}
 \theoremstyle{definition}
 
 \theoremstyle{remark}
 \newtheorem{rem}[theo]{Remark}
 \numberwithin{equation}{section}
 \newcommand{\apa}{\alpha}
 \newcommand{\bta}{\beta}
 \newcommand{\lma}{\lambda}
 \newcommand{\N}{\mathds{N}}
 \newcommand{\Z}{\mathds{Z}}
 \newcommand{\R}{\mathds{R}}
 \newcommand{\Cx}{\mathds{C}}
 \newcommand{\abar}{\bar{a}}
 \newcommand{\bbar}{\bar{b}}
 \newcommand{\spl}[1]{\displaystyle{#1}}
 \newcommand{\ups}[1]{\textsuperscript{#1}}
 \newcommand{\Iint}[2]{I_{#2}^{#1}}
 \newcommand{\cat}[2]{{}^C\!D_{#2}^{#1}}
 \newcommand{\catg}[2]{{}^C\!D_{#2}^{#1,\psi}}
 \newcommand{\der}[2]{{}^{RL}\!D_{#2}^{#1}}
 \newcommand{\derg}[2]{{}^{RL}\!D_{#2}^{#1,\psi}}
 \newcommand{\Erd}[2]{{}^{EK}\!D_{#2}^{#1}}
 \newcommand{\Had}[2]{{}^{CH}\!D_{#2}^{#1}}

\begin{document}

\title[On the fractional operators with respect to another function]
 {On the fractional operators with respect to another function}

\author{M. Benjemaa}

\address{%
Stability and Control of Systems,\\
University of Sciences of Sfax,\\
Tunisia}

\email{mondher.benjemaa@fss.usf.tn}


\author{F. Jerbi}
\address{Stability and Control of Systems,\\
University of Sciences of Sfax,\\
Tunisia}
\email{jerbiifatma@gmail.com}


\date{March 23, 2021}


\begin{abstract}
This paper is in concern with Cauchy problems involving the fractional derivatives with respect to another function. Results of existence, uniqueness, and Taylor series among others are established in appropriate functional spaces. We prove that these results are valid at once for several standard fractional operators such as the Riemann-Liouville and Caputo operators, the Hadamard operators, the Erd\'elyi-Kober operators, etc., depending on the choice of the scaling function. We also show that our technique can be useful to solve a wide range of Volterra integral equations. The numerical approximation of solutions of systems involving the fractional derivatives with respect to another function is also investigated and the optimal convergence rate of the schemes is reached in graded meshes, even in the case of singular solutions. Various examples and numerical tests, with an application to the Erd\'elyi-Kober operators, are performed at the end to illustrate the efficiency of proposed approach.
\end{abstract}

\maketitle

\section{Introduction}
The fractional derivatives provide an alternative framework to model various physical phenomena which involve nonlocal features \cite{Cap71,Machado17,Podlubny99}, such as visco-elasticity, seismology, chemistry, control theory, engineering, etc. \cite{Cap92,Fukunaga15,Old74}. Among the fractional operators one can cite the derivative in the sense of Riemann-Liouville (RL), the Caputo derivative, the Erd\'elyi-Kober derivative, the Hadamard derivative, etc. (see \cite{Diethelm10,Herrmann14,Hilfer00,Kilbas06,Kiryakova94,Lak09,Mil93,Podlubny99,Samko93,Sne79} among many others and references therein).\\

In view of the multitude of the fractional operator's definitions, and especially of some emergent ones, it is important to figure out if any connection between them exists, since each fractional operator is generally treated separately from the others, with the use of new definitions and properties, and in some cases, with similar or slightly modified existing proofs. Hence, finding a unifying framework to deal with many apparently different but possibly related operators is of interest to avoid such a redundancy. The fractional derivatives and integrals with respect to another function turn out to be one of these general classes of fractional operators. Introduced by Erd\'elyi in the sixties \cite{Erd64,Erd65}, they were extensively studied by Thomas J. Osler in the seventies in a series of papers \cite{Osl70a,Osl70b,Osl71,Osl72,Osl73}. Such operators, as well as their recent extension the $\psi$-Hilfer operators \cite{Hilfer00}, are regaining more interest nowadays \cite{Agr12,Alm17,Alm19,Sou18,Yan20} since they represent a generalization of several classical fractional derivatives, including the Riemann-Liouville derivative, the Hadamard derivative, the so-called generalized fractional derivative \cite{Benjemaa18,Katugampola14,Zen17} among others (we refer the reader to \cite{Sou18} for a more complete list). From a physical point of view, the concept of fractional derivatives with respect to another function have successfully been used to derive a generalization of the Scott-Blair models with time varying viscosity \cite{Col18}. In \cite{Gar18}, a new approach to the fractional Dodson diffusion equation using the fractional derivative with respect to another function is considered in order to get a deeper understanding of the memory effects in complex diffusion phenomena. In \cite{Alm18}, the authors showed that a system involving fractional derivative with respect to another function is more suitable to model the GDP\footnote{The gross domestic product} growth rate in the USA. We may also cite the work of O.P. Agrawal who introduced in \cite{Agr12} a fractional derivative with respect to two functions, called weight and scale functions, and several models using such an operator have been investigated, including the fractional diffusion equation \cite{Yufen13}, the fractional advection-diffusion model to describe the transport of a solute in aquifers \cite{Yuf13}, the generalized form of the Burger equation which can arise in several natural processes such as traffic flow, gas dynamics modeling, etc. \cite{Yufe13}. The study of the fractional operators with respect to other functions might thus be helpful to better model, understand, and unify the properties of several fractional operators, depending on the choice of the scaling function.\\

In this work, we propose to extensively investigate the fractional integrals and derivatives with respect to another function. First, we introduce an adequate functional framework in which these types fractional operators are well defined. Next, we derive some properties related to these operators, and we establish several results of existence and uniqueness of the solutions of Cauchy problems involving such a fractional derivatives. The cases of derivatives in the sense of Riemann-Liouville and Caputo are treated separately. These results allow us to make benefits of existing numerical schemes to approach the fractional integrals and derivatives with respect to a given function $\psi$, while keeping the optimal convergence rate, i.e. independently of the scaling function $\psi$. As an application, we were able to accurately approximate systems involving the Erd\'elyi-Kober operators, where the solutions or their first derivatives might be singular at the lower terminal point.\\

The paper consists of six sections. Some preliminaries are given in section \ref{SectionPre}. Next, we prove in section \ref{SectionRL} several properties in concern with the RL fractional operators with respect to another function, and various existing and new results are  extended or proven in suitable functional spaces. Similar results for the fractional derivatives with respect to another function in the sense of Caputo are derived in section \ref{SectionCap}. In section \ref{SectionErd}, we show that our approach can be applied to the Erd\'elyi-Kober's operators, and several properties for such a fractional operators are established. We also give some examples where Volterra integral equations of the first and the second kind are solved. Finally, we perform in section \ref{SectionNum} several numerical tests and we show that optimal rates of convergence can always be reached on graded meshes, even in case of singular solutions.

\section{Preliminaries}\label{SectionPre}
Let $b>a\geq 0$ and $\apa\geq 0$. The Riemann-Liouville (RL) fractional integral of order $\apa$ is defined for a function $u\in L^1(a,b)$ by
\begin{align}\label{LeftIntRL}
\Iint{\apa}{a}u(t) = \dfrac{1}{\Gamma(\apa)}\int_a^t \dfrac{u(s)}{(t-s)^{1-\apa}}\, ds.
\end{align}
Let $n\in \N$ (the set of positive integers) and let $AC^n[a,b]$ be the set of functions with an absolutely continuous $(n-1)$\ups{st} derivatives, i.e.
\begin{equation*}
AC^n[a,b] = \left\{u\in C^{n-1}[a,b],\ \exists\, \varphi \in L^1(a,b)\ / \ u^{(n-1)}(t) = u^{(n-1)}(a) + \int_a^t \varphi(s)\,ds\right\}.
\end{equation*}
Let $n=\lfloor \apa \rfloor + 1$. The (left) Riemann-Liouville fractional derivative of order $\apa$ and its Caputo modification are defined for any function $u\in AC^n[a,b]$ respectively by
\begin{equation}\label{LeftDerRL}
\der{\apa}{a}u(t) = \left(\Iint{n-\apa}{a}u\right)^{(n)}(t) = \dfrac{1}{\Gamma(n-\apa)} \left(\dfrac{d}{dt}\right)^n\int_a^t \dfrac{u(s)}{(t-s)^{\apa-n+1}}\ ds
\end{equation}
and
\begin{equation}\label{LeftDerCap}
\cat{\apa}{a}u(t) = \Iint{n-\apa}{a}\left(u^{(n)}\right)(t) = \dfrac{1}{\Gamma(n-\apa)} \int_a^t \dfrac{u^{(n)}(s)}{(t-s)^{\apa-n+1}}\ ds.
\end{equation}

More generally, if $\psi\in C^n[a,b]$ is a monotonously increasing function, then the Riemann-Liouville fractional integral and derivative of order $\apa$ with respect to the function $\psi$ are defined respectively by \cite{Kilbas06}
\begin{align}\label{LeftIntRL_g}
\Iint{\apa,\psi}{a}u(t) = \dfrac{1}{\Gamma(\apa)}\int_a^t \dfrac{\psi'(s)\,u(s)}{(\psi(t)-\psi(s))^{1-\apa}}\, ds
\end{align}
and
\begin{align}\label{LeftDerRL_g}
\derg{\apa}{a}u(t) & := \left(\dfrac{1}{\psi'(t)}\,\dfrac{d}{dt}\right)^n\left(\Iint{n-\apa,\psi}{a}u\right)(t) \nonumber \\
& \,= \dfrac{1}{\Gamma(n-\apa)} \left(\dfrac{1}{\psi'(t)}\,\dfrac{d}{dt}\right)^n\int_a^t \dfrac{\psi'(s)\,u(s)}{(\psi(t)-\psi(s))^{\apa-n+1}}\ ds.
\end{align}
In the sequel, we shall use the notations $\varrho := \frac{d}{d\psi} = \frac{1}{\psi'}\,\frac{d}{dt}$ and $\bar{\xi} = \psi(\xi)$ for simplicity. In \cite{Jar20}, the author proved that the Riemann-Liouville fractional derivatives with respect to $\psi$ are well defined on the set
\begin{equation*}
AC^n_\psi[a,b] := \left\{u:[a,b]\rightarrow \Cx\ \text{ s.t.} \ \varrho^{n-1}u\in AC[a,b]\right\}.
\end{equation*}
The Caputo fractional derivative operator with respect to another function is defined as \cite{Jar20}
\begin{align}\label{Cg_def}
\catg{\apa}{a}u(t) = \derg{\apa}{a}\left(u(t) - \sum_{k=0}^{n-1} \dfrac{\left(\varrho^{k}u\right)(a)}{k!}\big{(}\psi(t)-\psi(a)\big{)}^k\right).
\end{align}
In \cite{Alm17}, it is proven that for $u\in C^n[a,b]$, the Caputo fractional derivative operator with respect to another function can also be expressed as follows
\begin{align}\label{LeftDerCap_g}
\catg{\apa}{a}u(t) & = \Iint{n-\apa,\psi}{a}\left(\varrho^n u\right)(t) \nonumber \\
& \,= \dfrac{1}{\Gamma(n-\apa)} \int_a^t \dfrac{\psi'(s)\,(\varrho^n u)(s)}{(\psi(t)-\psi(s))^{\apa-n+1}}\ ds.
\end{align}
We shall see in the sequel that \eqref{LeftDerCap_g} holds true in the bigger space $AC^n_\psi[a,b]$ (see Appendix~A). Let us remark that if $\psi(t)=t$ (resp. $\psi(t) = \log t$, resp. $\psi(t) = \frac{t^\rho}{\rho}$ with $\rho>0$), then the Riemann-Liouville fractional derivative with respect to $\psi$ reduces to the standard Riemann-Liouville (resp. the Hadamard, resp. the so called generalized fractional derivative \cite{Katugampola14,Zen17}). The same claim is still true when considering the previous fractional derivatives in the Caputo sense.\\

Throughout all the paper, and unless specified, $\psi$ will denote a non-negative continuous function over an interval $(a,b]$, monotonously increasing and of class $C^n(a,b]$ such that $\psi'(t)\neq 0$ for all $t\in (a,b]$. Under these assumptions, $\psi$ is invertible and $\psi^{-1}\in C^n(\abar,\bbar]$, with the notation $\bar{\xi}=\psi(\xi)$. The following results will be useful in the sequel.
\begin{lemma}\label{lemmaDL}
Let $n\in \N$ and let $u$ and $\psi$ be a continuous functions, $n$ times derivable in $[a,b]$. Suppose $\psi$ is invertible and $\psi'(t)\neq 0$ for all $t\in [a,b]$. Then
\begin{align}\label{vn}
\varrho^n u = \left(u\circ \psi^{-1}\right)^{(n)}\circ \psi \quad \text{in } [a,b].
\end{align}
\end{lemma}
\begin{proof}
We prove the result by induction. The case $n=1$ is straightforward. Suppose \eqref{vn} is satisfied, then
\begin{align*}
\left(u\circ \psi^{-1}\right)^{(n+1)}\circ \psi & = \left(\left(\varrho^n u\right)\circ \psi^{-1}\right)'\circ \psi \\
   & = \left[\left(\dfrac{1}{\psi'}\left(\varrho^n u\right)'\right)\circ \psi^{-1}\right]\circ \psi \\
   & = \varrho^{n+1} u.
\end{align*}
\end{proof}

Let $1\leq p< \infty$ and define the space of $p$-integrable functions with respect to a function $\psi$:
\begin{equation}\label{Lpg}
L^p_\psi(a,b):=\left\{f:[a,b]\rightarrow \Cx,\ \int_a^b |f(s)|^p\,d\psi(s) < \infty\right\}.
\end{equation}

\begin{rem}\label{rmk_Lpg}
If $\psi'$ is bounded on $[a,b]$, then $L^p_\psi(a,b)=L^p(a,b)$.
\end{rem}

\begin{coro}\label{CoroEq}
Let $1\leq p< \infty$ and $n\in \N$. 
Then 
\begin{equation*}
u\in L^p_\psi(a,b) \Longleftrightarrow u\circ \psi^{-1}\in L^p(\abar,\bbar),
\end{equation*}
\begin{equation*}
u\in AC^n_\psi[a,b] \Longleftrightarrow u\circ \psi^{-1}\in AC^n[\abar,\bbar].
\end{equation*}
\end{coro}
\begin{proof}
The first assertion is immediate. Denote $z = u\circ \psi^{-1}$, then we obtain by Lemma~\ref{lemmaDL}
\begin{align*}
z\in AC^n[\abar,\bbar] & \Leftrightarrow \exists\, c\in \R \, \text{and}\, \varphi\in L^1(\abar,\bbar) \, \text{s.t.}\ z^{(n-1)}(\bar{t}) = c + \int_{\abar}^{\bar{t}} \varphi(s)\,ds \ \ \forall\ \bar{t}\in [\abar,\bbar] \\
& \Leftrightarrow z^{(n-1)}\circ \psi(t) = c + \int_{\psi(a)}^{\psi(t)} \varphi(s)\,ds \ \ \forall\ t\in [a,b] \\
& \Leftrightarrow \varrho^{n-1} u(t) = c + \int_{\psi(a)}^{\psi(t)} \varphi(s)\,ds \ \ \forall\ t\in [a,b] \\
& \Leftrightarrow \varrho^{n-1} u(t) = c + \int_{a}^{t} \psi'(s)\, \varphi\circ \psi(s)\,ds \ \ \forall\ t\in [a,b] \\
& \Leftrightarrow u \in AC^n_\psi[a,b].
\end{align*}
\end{proof}

\begin{propos}\label{prop_IDg}
Let $\apa\geq0$ and $n=\lceil \apa \rceil$. 
Then for any $u\in L^p_\psi(a,b)$, \mbox{$1\leq p < \infty$,} we have
\begin{enumerate}[leftmargin = 20pt, label = \roman*)]
\item $\Iint{\apa,\psi}{a}u = \left(\Iint{\apa}{\abar}\left(u\circ \psi^{-1}\right)\right)\circ \psi$,
\end{enumerate}
and for any $u\in AC^n_\psi[a,b]$,
\begin{enumerate}[leftmargin = 18pt]
\item[ii)] $\derg{\apa}{a}u = \left(\der{\apa}{\abar}\left(u\circ \psi^{-1}\right)\right)\circ \psi$,
\item[iii)] $\catg{\apa}{a}u = \left(\cat{\apa}{\abar}\left(u\circ \psi^{-1}\right)\right)\circ \psi.$
\end{enumerate}
\end{propos}
\begin{proof}
Let $u\in L^p_\psi(a,b)$. Then, using Corollary \ref{CoroEq}, we have $u\circ \psi^{-1}\in L^p(\abar,\bbar)$ and hence $\Iint{\apa}{\abar}\left(u\circ \psi^{-1}\right)$ 
is well defined on $[\abar,\bbar]$. It follows that for a.e. $t\in [a,b]$
\begin{align*}
\Iint{\apa,\psi}{a}u(t) & = \dfrac{1}{\Gamma(\apa)}\int_a^t \dfrac{\psi'(s)\,u(s)}{(\psi(t)-\psi(s))^{1-\apa}} ds = \dfrac{1}{\Gamma(\apa)}\int_{\abar}^{\bar{t}} \dfrac{u(\psi^{-1}(\xi))}{(\psi(t)-\xi)^{1-\apa}} d\xi\\ & = \Iint{\apa}{\abar}\left(u\circ \psi^{-1}\right)(\psi(t)).
\end{align*}
Now, let $u\in AC^n_\psi[a,b]$. Then, we have by Corollary \ref{CoroEq} that $u\circ \psi^{-1}\in AC^n[a,b]$ and hence $\der{\apa}{\abar}\left(u\circ \psi^{-1}\right)$ and $\cat{\apa}{\abar}\left(u\circ \psi^{-1}\right)$ are well defined on $[\abar,\bbar]$. It follows by $i)$ and Lemma \ref{lemmaDL} that for a.e. $t\in [a,b]$
\begin{align*}
\derg{\apa}{a}u(t) & = \varrho^n\left(\Iint{n-\apa,\psi}{a}u\right)(t) = \left(\dfrac{d}{dt}\right)^n\left[\left(\Iint{n-\apa,\psi}{a} u\right)\circ \psi^{-1}\right](\psi(t)) \\
& = \left(\dfrac{d}{dt}\right)^n\left[\Iint{n-\apa}{\abar}\left(u\circ \psi^{-1}\right)\right](\psi(t)) \\
& = \der{\apa}{\abar}\left(u\circ \psi^{-1}\right)(\psi(t)),
\end{align*}
and
\begin{align*}
\catg{\apa}{a}u(t) & = \Iint{n-\apa,\psi}{a}\left(\varrho^n u\right)(t) \\
& = \Iint{n-\apa}{\abar}\left(\varrho^n u\circ \psi^{-1}\right)(\psi(t)) \\
& = \Iint{n-\apa}{\abar}\left(u\circ \psi^{-1}\right)^{(n)}(\psi(t)) \\
& = \cat{\apa}{\abar}\left(u\circ \psi^{-1}\right)(\psi(t)).
\end{align*}
\end{proof}

\section{The Riemann-Liouville fractional operators with respect to another function}\label{SectionRL}
\begin{propos}
Let $p\geq 1$. Then, the fractional integral operator with respect to another function $I^{\apa,\psi}_a$ is bounded in $L^p_\psi(a,b)$:
\begin{align*}
\|I^{\apa,\psi}_a u\|_{L^p_\psi(a,b)} \leq K_\psi \|u\|_{L^p_\psi(a,b)} \quad \text{with } K_\psi = \dfrac{(\bbar-\abar)^\apa}{\Gamma(\apa+1)}.
\end{align*}
\end{propos}

\begin{proof}
First, remark that $\|f\circ \psi\|_{L^p_\psi(a,b)} = \|f\|_{L^p(\abar,\bbar)}$, with the notation $\bar{\xi}=\psi(\xi)$. Using Proposition \ref{prop_IDg} and the continuity of the operator $I^{\apa}_a$, $$ \|I^{\apa}_a f\|_{L^p(a,b)} \leq K \|f\|_{L^p(a,b)} \quad \text{with } K = \dfrac{(b-a)^\apa}{\Gamma(\apa+1)}$$ (see e.g. \cite[Lemma 2.1]{Kilbas06}), one obtain
\begin{align*}
\|I^{\apa,\psi}_a u\|_{L^p_\psi(a,b)} & = \left\|\left(\Iint{\apa}{\abar}\left(u\circ \psi^{-1}\right)\right)\circ \psi\right\|_{L^p_\psi(a,b)}\\
& = \left\|\Iint{\apa}{\abar}\left(u\circ \psi^{-1}\right)\right\|_{L^p(\abar,\bbar)}\\
& \leq K_\psi \left\|u\circ \psi^{-1}\right\|_{L^p(\abar,\bbar)}\\
& = K_\psi \left\|u\right\|_{L^p_\psi(a,b)}
\end{align*}
\end{proof}

\begin{propos}({\bf semi-group law})\label{prop_semigpe_RL}
Let $\apa > 0$ and $\bta > 0$. Then, for any $1\leq p < \infty$ and $u\in L^p_\psi(a,b)$
\begin{enumerate}
\item $\Iint{\apa,\psi}{a} \Iint{\bta,\psi}{a} u = \Iint{\apa+\bta,\psi}{a} u.$
\item $\derg{\apa}{a}\Iint{\apa,\psi}{a} u = u.$
\item $\derg{\bta}{a}\Iint{\apa,\psi}{a} u = \Iint{\apa-\bta,\psi}{a}u \ \ \forall \ \apa \geq \bta.$
\item Let $m\in \N$, then 
\begin{equation*}
\varrho^m\Iint{\apa,\psi}{a} u = \left\{
\begin{array}{ll}
\Iint{\apa-m,\psi}{a}u & \text{if }\, m \leq \apa \\
\derg{m-\apa}{a}u & \text{if }\, m \geq \apa.
\end{array}
\right.
\end{equation*}
\end{enumerate}
\end{propos}

\begin{proof}
We only prove the first assertion, since the proofs of the other identities follow a similar idea. Let $u\in L^p_\psi(a,b)$, then $u\circ \psi^{-1}\in L^p(\abar,\bbar)$. It follows from Lemma 2.3 in \cite{Kilbas06} and Proposition \ref{prop_IDg} that for a.e. $t\in [a,b]$
\begin{align*}
\Iint{\apa,\psi}{a} \Iint{\bta,\psi}{a} u(t) &= \Iint{\apa}{\abar} \left( \left(\Iint{\bta,\psi}{a} u\right)\circ \psi^{-1}\right)(\psi(t))\\
&= \Iint{\apa}{\abar} \left( \Iint{\bta}{\abar} (u\circ \psi^{-1})\right)(\psi(t))\\
&= \Iint{\apa +\bta}{\abar} ( u \circ \psi^{-1})(\psi(t))\\
&= \Iint{\apa +\bta,\psi}{a} u(t).
\end{align*}
\end{proof}

More generally, we have the following semi-group identity.
\begin{propos}\label{prop_permut_RL}
Let $\apa > 0$, $\bta > 0$ and $\lma$ and $\nu$ in $\R$. Then
\begin{equation*}
\Iint{\apa,\psi}{a} \left( \psi(t)^\lma\Iint{\bta,\psi}{a} \left(\psi(t)^\nu u(t)\right)\right) = \psi(t)^{\apa+\lma}\Iint{\bta,\psi}{a} \left( \psi(t)^{-(\apa+\bta+\lma)}\Iint{\apa,\psi}{a} \left(\psi(t)^{\bta+\lma+\nu}u(t)\right)\right).
\end{equation*}
\end{propos}
\begin{proof}
First, we rewrite $\Iint{\apa,\psi}{a}u$ as follows:
\begin{equation*}
\Iint{\apa,\psi}{a}u(t) = \dfrac{\psi(t)^\apa}{\Gamma(\apa)}\int_{\frac{\psi(a)}{\psi(t)}}^1 (1-x)^{\apa-1} \, u\circ \psi^{-1}(x\psi(t))\,dx,
\end{equation*}
which can be obtained by setting $x = \psi(s)/\psi(t)$ in \eqref{LeftIntRL_g}. It follows that
\begin{equation*}
\Iint{\bta,\psi}{a} \left(\psi^\nu\, u\right)(\xi) = \dfrac{\psi(\xi)^{\bta+\nu}}{\Gamma(\bta)}\int_{\frac{\psi(a)}{\psi(\xi)}}^1 (1-x)^{\bta-1}\,x^\nu \, u\circ \psi^{-1}(x\psi(\xi))\,dx.
\end{equation*}
Using the Fubini's formula, we get
\begin{align*}
& \Iint{\apa,\psi}{a} \left( \psi(t)^\lma\Iint{\bta,\psi}{a} \left(\psi(t)^\nu u(t)\right)\right) \\
& \quad = \dfrac{\psi(t)^\apa}{\Gamma(\apa)}\int_{\frac{\psi(a)}{\psi(t)}}^1 (1-x)^{\apa-1} \, (x\psi(t))^\lma\,\Iint{\bta,\psi}{a} \left(\psi^\nu u\right)\left(\psi^{-1}(x\psi(t))\right) \,dx \\
& \quad = \dfrac{\psi(t)^{\apa+\lma}}{\Gamma(\apa)}\int_{\frac{\psi(a)}{\psi(t)}}^1 (1-x)^{\apa-1} \,x^\lma\,\dfrac{(x\psi(t))^{\bta+\nu}}{\Gamma(\bta)}\int_{\frac{\psi(a)}{x\psi(t)}}^1 (1-y)^{\bta-1}\,y^\nu \,u\circ \psi^{-1}(yx\psi(t))\,dy\, dx \\
& \quad = \dfrac{\psi(t)^{\apa+\bta+\lma+\nu}}{\Gamma(\apa)\Gamma(\bta)}\int_{\frac{\psi(a)}{\psi(t)}}^1 (1-x)^{\apa-1} \,x^{\bta+\lma+\nu}\,\int_{\frac{\psi(a)}{x\psi(t)}}^1 (1-y)^{\bta-1}\,y^\nu \,u\circ \psi^{-1}(yx\psi(t))\,dy\, dx \\
& \quad = \dfrac{\psi(t)^{\apa+\bta+\lma+\nu}}{\Gamma(\apa)\Gamma(\bta)}\int_{\frac{\psi(a)}{\psi(t)}}^1 (1-y)^{\bta-1}\,y^\nu \,\int_{\frac{\psi(a)}{y\psi(t)}}^1 (1-x)^{\apa-1} \,x^{\bta+\lma+\nu}\,u\circ \psi^{-1}(yx\psi(t))\,dx\, dy \\
& \quad = \psi(t)^{\apa+\lma}\Iint{\bta,\psi}{a} \left( \psi(t)^{-(\apa+\bta+\lma)}\Iint{\apa,\psi}{a} \left(\psi(t)^{\bta+\lma+\nu}u(t)\right)\right).
\end{align*}
\end{proof}

\begin{propos}\label{prop_I-beta}
Let $\apa>0$ with $\apa\not\in \N$ and $n=\lceil \apa \rceil$. Let $\psi$ be a monotonously increasing function of class $C^n[a,b]$ such that $\psi'(t)\neq 0$ for all $t\in [a,b]$. Then for any $u\in C^n[a,b]$ and any $t\in [a,b]$
\begin{equation*}
\derg{\apa}{a} u(t) = \Iint{-\apa,\psi}{a} u(t).
\end{equation*}
\end{propos}

\begin{proof}
From the hypothesis we have $u\circ \psi^{-1}\in C^n[\abar,\bbar]$. Let $t\in [a,b]$, then using Lemma 2.21 in \cite{Diethelm10} and Proposition \ref{prop_IDg}, we obtain
\begin{align*}
\derg{\apa}{a} u(t) & = \der{\apa}{\psi(a)} \left(u\circ \psi^{-1}\right)(\psi(t))\\
& = \dfrac{1}{\Gamma(-\apa)}\int_{\psi(a)}^{\psi(t)} \left(\psi(t)-s\right)^{-\apa-1}\left(u\circ \psi^{-1}\right)(s)\,ds\\
& = \dfrac{1}{\Gamma(-\apa)}\int_{a}^t \left(\psi(t)-\psi(\xi)\right)^{-\apa-1}\,u(\xi)\,\psi'(\xi)\,d\xi\\
& = \Iint{-\apa,\psi}{a} u(t).
\end{align*}
\end{proof}

\begin{propos}\label{prop_IagDag}
Let $\apa > 0$ and $n=\lceil \apa \rceil$. If $u\in L^1_\psi(a,b)$ and $\Iint{n-\apa,\psi}{a} u \in AC^n_\psi[a,b]$, then for a.e. $t\in [a,b]$
\begin{align*}
\Iint{\apa,\psi}{a}\left( \derg{\apa}{a}u\right)(t) = u(t) - \sum_{k=1}^n \dfrac{\derg{\apa-k}{a}u(a)}{\Gamma(\apa-k+1)}(\psi(t)-\psi(a))^{\apa-k}.
\end{align*}
\end{propos}
\begin{proof}
By the hypothesis and Corollary \ref{CoroEq}, we have $u\circ \psi^{-1}\in L^1(\abar,\bbar)$ and $\left(\Iint{n-\apa,\psi}{a} u\right)\circ \psi^{-1} \in AC^n[\abar,\bbar]$. From Proposition \ref{prop_IDg} we deduce $\Iint{n-\apa}{\abar}\left( u\circ \psi^{-1}\right) \in AC^n[\abar,\bbar]$. Using Proposition \ref{prop_IDg} and Lemma \ref{lemmaDL} hereinabove and Lemma~2.5 in \cite{Kilbas06}, we obtain
\begin{align*}
\Iint{\apa,\psi}{a}\left( \derg{\apa}{a}u\right)(t) & = \Iint{\apa}{\abar}\left[\left(\derg{\apa}{a}u\right)\circ \psi^{-1}\right](\psi(t))\\
& = \Iint{\apa}{\abar}\left[\der{\apa}{\abar}\left(u\circ \psi^{-1}\right)\right](\psi(t))\\
& = u\circ \psi^{-1}(\psi(t)) - \sum_{k=1}^n \dfrac{\left[\Iint{n-\apa}{\abar}\left(u\circ \psi^{-1}\right)\right]^{(n-k)}(\abar)}{\Gamma(\apa-k+1)}(\psi(t)-\abar)^{\apa-k} \\
& = u(t) - \sum_{k=1}^n \dfrac{\left[\left(\Iint{n-\apa,\psi}{a}u\right)\circ \psi^{-1}\right]^{(n-k)}(\abar)}{\Gamma(\apa-k+1)}(\psi(t)-\psi(a))^{\apa-k} \\
& = u(t) - \sum_{k=1}^n \dfrac{\varrho^{n-k}\left(\Iint{n-\apa,\psi}{a}u\right)(a)}{\Gamma(\apa-k+1)}(\psi(t)-\psi(a))^{\apa-k}\\
& = u(t) - \sum_{k=1}^n \dfrac{\derg{\apa-k}{a}u(a)}{\Gamma(\apa-k+1)}(\psi(t)-\psi(a))^{\apa-k},
\end{align*}
where the last equality is obtained using Proposition \ref{prop_semigpe_RL}.
\end{proof}

\begin{propos}\label{prop_Da_Db}
Let $\apa>0$ and $\bta>0$ such that $n-1< \apa \leq n$ and $m-1< \bta \leq m$ with $n, \, m\in \N$. Let $u\in L^1_\psi(a,b)$ such that $\Iint{m-\bta,\psi}{a} u \in AC^m_\psi[a,b]$, then for a.e. $t\in [a,b]$
\begin{align*}
\derg{\apa}{a}\left(\derg{\bta}{a}u\right)(t) & = \derg{\apa+\bta}{a}u(t) - \sum_{k=1}^{m-1} \dfrac{\derg{\bta-k}{a}u(a)}{\Gamma(1-k-\apa)}\left(\psi(t)-\psi(a)\right)^{-k-\apa} \\
& \qquad - \dfrac{(\psi(t)-\psi(a))^{-m-\apa}}{\Gamma(1-m-\apa)}\,\lim_{t\to a}\Iint{m-\bta,\psi}{a}u(t).
\end{align*}
\end{propos}

\begin{proof}
The proof is similar to that of Proposition \ref{prop_IagDag} and is essentially based on Proposition \ref{prop_IDg}, Lemma \ref{lemmaDL} hereinabove and Property 2.4 in \cite{Kilbas06}. 
\end{proof}

\begin{propos}({\bf Leibniz formula})\label{prop_Leibniz}
Let $\apa \in \R_+\backslash \N_0$, and assume that $u$, $v$ and $\psi$ along all their derivatives are continuous on $[a,b]$, with $\psi'(t)\neq 0$ for all $t\in [a,b]$. Then, for any $t\in [a,b]$
\begin{equation*}
\derg{\apa}{a}(uv)(t) = \sum_{k=0}^{\infty} \begin{pmatrix}\apa\\k\end{pmatrix} \varrho^k u(t)\, \derg{\apa-k}{a}v(t),
\end{equation*}
where $\begin{pmatrix}\apa\\k\end{pmatrix} = \dfrac{\Gamma(\apa+1)}{\Gamma(k+1)\,\Gamma(\apa-k+1)}$.
\end{propos}
\begin{proof}
A direct consequence of the Leibniz formula for the standard RL fractional derivative (see e.g. \cite[eq (2.202)]{Podlubny99}), Proposition \ref{prop_IDg} and Lemma~\ref{lemmaDL}.
\end{proof}

More generally, we have the following (symmetric) product rule.
\begin{propos}
Assume the hypothesis of Proposition \ref{prop_Leibniz}, and let $\nu$ be an arbitrary real (or complex) number, $\nu\not \in \Z_0$. Then, for any $t\in [a,b]$
\begin{equation*}
\derg{\apa}{a}(uv)(t) = \sum_{k=-\infty}^{\infty} \begin{pmatrix}\apa\\k+\nu\end{pmatrix} \derg{k+\nu}{a} u(t)\, \derg{\apa-k-\nu}{a}v(t).
\end{equation*}
\end{propos}
\begin{proof}
The proof follows from the product rule for the standard RL fractional derivative (see e.g. \cite[Theorem 1]{Osl70b}), Proposition \ref{prop_IDg} and Lemma \ref{lemmaDL}.
\end{proof}


\begin{lemma}\label{lem_der}
Let $\apa > 0$ and $n = \lceil \apa \rceil$. Then, for $t>a$
\begin{align*}
& \Iint{\apa,\psi}{a}\left(\psi(t)-\psi(a)\right)^\bta = \dfrac{\Gamma(1+\bta)}{\Gamma(1+\bta+\apa)}\left(\psi(t)-\psi(a)\right)^{\bta+\apa}, & \bta > -1,\\
& \derg{\apa}{a}\left(\psi(t)-\psi(a)\right)^\bta = \dfrac{\Gamma(1+\bta)}{\Gamma(1+\bta-\apa)}\left(\psi(t)-\psi(a)\right)^{\bta-\apa}, & \bta > -1,\\
& \derg{\apa}{a}\left(\psi(t)-\psi(a)\right)^{\apa-i} = 0, & i=1,\dots n.
\end{align*}
\end{lemma}
\begin{proof}
A direct consequence of Proposition \ref{prop_IDg} hereinabove and Property~2.1 in \cite{Kilbas06}.
\end{proof}

For power functions that do not depend on the lower terminal, we prove the following:
\begin{lemma}\label{lem_der_power}
Let $\apa > 0$ and $n = \lceil \apa \rceil$. If $\psi(a)\neq 0$, then for $\bta>-1$ and $t>a$
\begin{align*}
& \Iint{\apa,\psi}{a} \psi(t)^\bta = \sum_{k\geq 0}\dfrac{\Gamma(1+\bta)}{\Gamma(\bta-k)\Gamma(1+k+\apa)}\psi(a)^{\bta-k}\left(\psi(t)-\psi(a)\right)^{k+\apa},\\
& \derg{\apa}{a} \psi(t)^\bta = \sum_{k\geq 0}\dfrac{\Gamma(1+\bta)}{\Gamma(\bta-k)\Gamma(1+k-\apa)}\psi(a)^{\bta-k}\left(\psi(t)-\psi(a)\right)^{k-\apa}.
\end{align*}
In particular, for $m\in \N_0$
\begin{align*}
& \Iint{\apa,\psi}{a} \psi(t)^m = \sum_{k=0}^m\dfrac{m!}{(m-k)!\,\Gamma(1+k+\apa)}\psi(a)^{m-k}\left(\psi(t)-\psi(a)\right)^{k+\apa},\\
& \derg{\apa}{a} \psi(t)^m = \sum_{k=0}^m\dfrac{m!}{(m-k)!\,\Gamma(1+k-\apa)}\psi(a)^{m-k}\left(\psi(t)-\psi(a)\right)^{k-\apa}.
\end{align*}
\end{lemma}
\begin{proof}
Let $b\neq 0$, then a Taylor expansion of the function $x\mapsto x^\bta$ writes
$$x^\bta = \sum_{k\geq 0} \dfrac{\Gamma(1+\bta)}{k!\,\Gamma(\bta-k)}b^{\bta-k} (x-b)^k.$$
Substitute $x=\psi(t)$, $b=\psi(a)$ and use Lemma \ref{lem_der} achieves the proof.
\end{proof}

Let $\apa>0$ and $n= \lceil \apa \rceil$. Consider the following fractional differential system
\begin{align}\label{fracds_RL}
\left\{
\begin{array}{l}
\derg{\apa}{a}u(t) = f(t,u(t)), \quad t\in [a,b] \\ \\
\derg{\apa-k}{a}u(a) = a_k, \ \ k= 1,\dots ,\, n-1, \ \ \spl{\lim_{t\to a}}\Iint{n-\apa,\psi}{a}u(t) = a_n.
\end{array}
\right.
\end{align}
with $f$ a given function and $a_k\in \R$ for $k=1,\dots ,\, n$. We have the following integral representation of the solution of \eqref{fracds_RL}.
\begin{propos}\label{prop_volt}
Let $G\subseteq \R$ be an open set and assume $f : (a,b]\times G \rightarrow \R$ is a function such that $t\mapsto f(t,\cdot)\in L^1_\psi(a,b)$. Then a function $u\in L^1_\psi(a,b)$ is a solution of \eqref{fracds_RL} if and only if $u$ is a solution of the non-linear second kind Volterra integral equation
\begin{align}\label{u=volterra_RL}
u(t) = \sum_{k=1}^{n}\dfrac{a_k}{\Gamma(\apa-k+1)}\big{(}\psi(t)-\psi(a)\big{)}^{\apa-k}  + \dfrac{1}{\Gamma(\apa)}\int_a^t \dfrac{\psi'(s)\,f(s,u(s))}{\left(\psi(t)-\psi(s)\right)^{1-\apa}}\,ds
\end{align}
with $a_k = \derg{\apa-k}{a}u(a)$ for $k=1,\dots ,\, n-1$, and $a_n=\spl{\lim_{t\to a}}\Iint{n-\apa,\psi}{a}u(t)$.
\end{propos}
\begin{proof}
Using Proposition \ref{prop_IDg}, one can show that $u$ is a solution of \eqref{fracds_RL} if and only if $v = u\circ \psi^{-1}$ is a solution of the system
\begin{align}\label{fracds_RL_v}
\left\{
\begin{array}{l}
\der{\apa}{\abar}v(x) = F(x,v(x)), \quad x\in [\abar,\bbar] \\ \\
\der{\apa-k}{\abar}v(\abar) = a_k, \ \ k=1,\dots ,\, n-1, \ \ \spl{\lim_{x\to \abar}}\Iint{n-\apa}{\abar}v(x) = a_n.
\end{array}
\right.
\end{align}
with $F(x,y) = f(\psi^{-1}(x),y)$. Noticing that $t\mapsto f(t,\cdot)\in L^1_\psi(a,b) \Leftrightarrow x\mapsto F(x,\cdot) \in L^1(\abar,\bbar)$ and $u\in L^1_\psi(a,b)\Leftrightarrow v = u\circ \psi^{-1} \in L^1(\abar,\bbar)$, and using Theorem 3.1 in \cite{Kilbas06}, we deduce that $v$ is a solution of \eqref{fracds_RL_v} if and only if $v$ satisfies for a.e. $x\in [\abar,\bbar]$
\begin{equation}\label{v=(x)}
v(x) =  \sum_{k=1}^{n}\dfrac{a_k}{\Gamma(\apa-k+1)}\big{(}x-\abar\big{)}^{\apa-k}  + \dfrac{1}{\Gamma(\apa)}\int_{\abar}^x \dfrac{F(s,v(s))}{\left(x-s\right)^{1-\apa}}\,ds.
\end{equation}
Finally, the result follows by taking $x=\psi(t)$ in equation \eqref{v=(x)}.
\end{proof}
Define the space
\begin{equation}\label{Lalphag}
L^\apa_\psi(a,b) := \left\{\varphi\in L^1_\psi(a,b),\ \derg{\apa}{a}\varphi\in L^1_\psi(a,b)\right\}
\end{equation}
where $L^1_\psi(a,b)$ is given in \eqref{Lpg}. Then we have the following result.
\begin{theo}\label{exist_uniq_RL}
Let $\apa>0$ and $n=\lceil\apa\rceil$. Let $G\subseteq \R$ be an open set and let $f : (a,b]\times G \rightarrow \R$ be a function such that $t\mapsto f(t,\cdot)\in L^1_\psi(a,b)$. Assume that $f$ fulfills a Lipschitz condition with respect to its second variable. Then the Cauchy problem \eqref{fracds_RL} admits a unique solution $u\in L^\apa_\psi(a,b)$.
\end{theo}
\begin{proof}
We have established in Proposition \ref{prop_volt} that $u$ is a solution of \eqref{fracds_RL} if and only if $v:=u\circ \psi^{-1}$ is a solution of \eqref{fracds_RL_v}. Since $(t,y)\mapsto F(t,y):=f(\psi^{-1}(t),y)$ is Lipschitzian with respect to its second variable, we obtain from Theorem 3.3 in \cite{Kilbas06} that the system \eqref{fracds_RL_v} admits a unique solution $v\in L^\apa(\abar,\bbar):= \left\{\varphi\in L^1(\abar,\bbar),\ \der{\apa}{\abar}\varphi\in L^1(\abar,\bbar)\right\}$. Finally, the result follows by noticing that $v\in L^\apa(\abar,\bbar) \Leftrightarrow u = v\circ \psi \in L^\apa_\psi(a,b)$.
\end{proof}

\begin{coro}\label{CoroEq_RL}
A function $u\in L^1_\psi(a,b)$ is a solution of \eqref{fracds_RL} if and only if $u = v \circ \psi$ with $v\in L^1(\abar,\bbar)$ is a solution of the Riemann-Liouville fractional differential system
\begin{align*}
\left\{
\begin{array}{l}
\der{\apa}{\abar}v(t) = F(t,v(t)), \quad t\in [\abar,\bbar] \\ \\
\der{\apa-k}{\abar}v(\abar) = a_k, \quad k=1,\dots ,\, n-1, \ \ \spl{\lim_{t\to \abar}}\Iint{n-\apa}{\abar}v(t) = a_n.
\end{array}
\right.
\end{align*}
with the notation $\bar{\xi} = \psi(\xi)$ and $F(t,x) = f(\psi^{-1}(t),x)$.
\end{coro}
\begin{rem}
Corollary \ref{CoroEq_RL} allows us to easily derive high order numerical schemes for many fractional operators, such as the Hadamard fractional operators \cite{Kilbas06}, the generalized operators \cite{Katugampola14,Zen17}, the Erd\'elyi-Kober fractional operators \cite{Sne75}, etc. (see Section \ref{SectionNum} for more details)
\end{rem}

\section{The Caputo fractional operators with respect to another function}\label{SectionCap}
\begin{propos}\label{prop_Cg=Igamma}
Let $\apa > 0$, $n=\lceil \apa \rceil$, and let $u\in AC^n_\psi[a,b]$. Then for a.e. $t\in [a,b]$
\begin{align*}
\catg{\apa}{a}u(t) & = \Iint{n-\apa,\psi}{a}\left(\varrho^n u\right)(t).
\end{align*}
In particular, $\catg{n}{a}u = \varrho^n u$ for any $n\in \N$.
\end{propos}
\begin{proof}
Using Corollary \ref{CoroEq}, we have $u\circ \psi^{-1}\in AC^n[\abar,\bbar]$. It follows from \cite[Theorem 2.1]{Kilbas06}, Proposition \ref{prop_IDg} and Lemma \ref{lemmaDL} that for a.e. $t\in [a,b]$
\begin{align*}
\catg{\apa}{a}u(t) & = \cat{\apa}{\abar}\left(u\circ \psi^{-1}\right)(\psi(t)) \\
& = \Iint{n-\apa}{\abar}\left(\left(u\circ \psi^{-1}\right)^{(n)}\right)(\psi(t)) \\
& = \Iint{n-\apa}{\abar}\left(\left(\varrho^n u\right)\circ \psi^{-1}\right)(\psi(t)) \\
& = \Iint{n-\apa,\psi}{a}\left(\varrho^n u\right)(t).
\end{align*}
\end{proof}
\begin{propos}({\bf Composition})\label{PropCompo}
Let $\apa>0$ and $n= \lceil \apa \rceil$. Then, for $u\in C[a,b]$
\begin{align*}
\catg{\apa}{a}\left(\Iint{\apa,\psi}{a}u\right)(t) = u(t),
\end{align*}
and for $u\in AC^n_\psi[a,b]$ 
\begin{align*}
\Iint{\apa,\psi}{a}\left(\catg{\apa}{a}u\right)(t) = u(t) - \sum_{k=0}^{n-1} \dfrac{\left(\varrho^k u\right)(a)}{k!}\big{(}\psi(t)-\psi(a)\big{)}^k.
\end{align*}
\end{propos}
\begin{proof}
Since $\psi$ is continuous, then $u\circ \psi^{-1}$ is also continuous. The rest of the proof is an immediate consequence of Theorems 3.7 and 3.8 in \cite{Diethelm10} and Proposition \ref{prop_IDg}, Lemma \ref{lemmaDL} and Corollary \ref{CoroEq} hereinabove. 
\end{proof}
\begin{theo}\label{Theo_continu_CDa}
Let $\apa > 0$, $n=\lceil \apa \rceil$, and let $u$ and $\psi$ in $C^n[a,b]$ such that $\psi'(t)\neq 0$ for all $t\in [a,b]$. Then $\catg{\apa}{a}u\in C[a,b]$. In particular, if $\apa\not \in \N_0$ then $\catg{\apa}{a}u(a)=0$. 
\end{theo}
\begin{proof}
Since $\psi^{-1}\in C^n[\abar,\bbar]$ then $u\circ \psi^{-1}\in C^n[\abar,\bbar]$. Using that $C^n[a,b]\subset AC^n_\psi[a,b]$ (see Appendix \ref{SectionAppA}), we obtain by Theorem 2.2 in \cite{Kilbas06} and Proposition~\ref{prop_IDg} that $\catg{\apa}{a}u = \left(\cat{\apa}{\abar}\left(u\circ \psi^{-1}\right)\right)\circ \psi \in C[a,b]$. Moreover, we have from Theorem~2.2 in \cite{Kilbas06} that $\cat{\apa}{\abar}\left(u\circ \psi^{-1}\right)(\abar) = 0$ if $\apa\not \in \N_0$. Applying again Proposition~\ref{prop_IDg} yields $\catg{\apa}{a}u(a)=0$.
\end{proof}
\begin{rem}
Theorem \ref{Theo_continu_CDa} may fails if one suppose only $\psi'\neq 0$ in $(a,b]$ (instead of $\psi'\neq 0$ in $[a,b]$). As a counterexample, one may consider $u(t) = t$ and $\psi(t)=t^{1/\apa}$ with $\apa\in (0,1)$, yielding $\catg{\apa}{0}u(t) = \Gamma(1+\apa)\ \forall\, t\geq 0$ (see Lemma \ref{lem_der_C} below), and hence $\catg{\apa}{0}u(0)\neq 0$.
\end{rem}
\begin{theo}({\bf Taylor's series})
Let $\apa \in (0,1]$ and let $m$ be an arbitrary non-negative integer. Suppose $\catg{k\apa}{a}u\in C[a,b]$ for $k=0,1\dots m+1$, where $\catg{k\apa}{a} := \catg{\apa}{a}\circ\dots\circ\catg{\apa}{a}$ ($k$ times) is the $k$\ups{th} sequential fractional derivative operator. Then 
\begin{align*}
u(t) = \sum_{j=0}^m \dfrac{\catg{j\apa}{a}u(a)}{\Gamma(j\apa+1)} \left(\psi(t)-\psi(a)\right)^{j\apa} + \dfrac{\catg{(m+1)\apa}{a}u(\xi)}{\Gamma((m+1)\apa+1)}\left(\psi(t)-\psi(a)\right)^{(m+1)\apa}
\end{align*} 
with $\xi \in\, (a,t)$.
\end{theo}
\begin{proof}
The result can be derived by using Theorem 3 in \cite{Odi07} and Proposition~\ref{prop_IDg} hereinabove.
\end{proof}
\begin{lemma}\label{lem_der_C}
Let $\apa > 0$ and $n = \lceil \apa\rceil$. Then
\begin{align*}
&\catg{\apa}{a}\big{(}\psi(t)-\psi(a)\big{)}^\bta = \dfrac{\Gamma(1+\bta)}{\Gamma(1+\bta-\apa)}\big{(}\psi(t)-\psi(a)\big{)}^{\bta-\apa},& \bta > n-1,\\
&\catg{\apa}{a}\big{(}\psi(t)-\psi(a)\big{)}^k = 0, & k=0,\, 1,\dots,\, n-1.
\end{align*}
\end{lemma}
\begin{proof}
A direct consequence of Proposition \ref{prop_IDg} hereinabove and Property~2.16 in \cite{Kilbas06}.
\end{proof}

\begin{lemma}\label{lem_der_C_power}
Let $\apa > 0$ and $n = \lceil \apa \rceil$. If $\psi(a)\neq 0$, then for $\bta>n-1$ and $t>a$
\begin{align*}
\catg{\apa}{a} \psi(t)^\bta = \sum_{k\geq n}\dfrac{\Gamma(1+\bta)}{\Gamma(\bta-k)\Gamma(1+k-\apa)}\psi(a)^{\bta-k}\left(\psi(t)-\psi(a)\right)^{k-\apa}.
\end{align*}
In particular, for any integer $m\geq n$
\begin{align*}
\catg{\apa}{a} \psi(t)^m = \sum_{k=n}^m\dfrac{m!}{(m-k)!\,\Gamma(1+k-\apa)}\psi(a)^{m-k}\left(\psi(t)-\psi(a)\right)^{k-\apa}.
\end{align*}
\end{lemma}
\begin{proof}
The proof is omitted since it is similar to that of Lemma \ref{lem_der_power}.
\end{proof}

Let $\apa>0$ and $n= \lceil \apa \rceil$. Consider the following Caputo fractional differential system with respect to another function
\begin{align}\label{fracds_c}
\left\{
\begin{array}{l}
\catg{\apa}{a}u(t) = f(t,u(t)), \quad t\in [a,b] \\ \\
\left(\varrho^k u\right)(a) = a_k, \quad k=0,\, 1,\dots ,\, n-1
\end{array}
\right.
\end{align}
with $f$ a given function and $a_k\in \R$ for $k=0,\, 1,\dots ,\, n-1$.
\begin{theo}\label{TheoEqCap}
A function $u\in AC^n_\psi[a,b]$ is a solution of \eqref{fracds_c} if and only if $u = v \circ \psi$ with $v\in AC^n[\abar,\bbar]$ is a solution of the Caputo differential system
\begin{align}\label{fde1}
\left\{
\begin{array}{l}
\cat{\apa}{\abar}v(t) = F(t,v(t)), \quad t\in [\abar,\bbar] \\ \\
v^{(k)}(\abar) = a_k, \quad k=0,\, 1,\dots ,\, n-1
\end{array}
\right.
\end{align}
with $F(t,x) = f(\psi^{-1}(t),x)$.
\end{theo}
\begin{proof}
A direct consequence of Proposition \ref{prop_IDg}, Lemma \ref{lemmaDL} and Corollary \ref{CoroEq}.
\end{proof}
\begin{propos}\label{prop_volt_c}
Assume $f$ is continuous over $[a,b]\times\R$. Then a function $u\in AC^n_\psi[a,b]$ is a solution of \eqref{fracds_c} if and only if $u$ is a solution of the of the non-linear second kind Volterra integral equation
\begin{align}\label{u=volterra_C}
u(t) = \sum_{k=0}^{n-1}\dfrac{a_k}{k!}\big{(}\psi(t)-\psi(a)\big{)}^k  + \dfrac{1}{\Gamma(\apa)}\int_a^t \dfrac{\psi'(s)\,f(s,u(s))}{\left(\psi(t)-\psi(s)\right)^{1-\apa}}\,ds
\end{align}
with $a_k = \left(\varrho^k u\right)(a)$ for $k=0,\, 1,\dots ,\, n-1$.
\end{propos}
\begin{proof}
"$\Rightarrow$" Apply the fractional integral with respect to another function \eqref{LeftIntRL_g} on both sides of equation \eqref{fracds_c} and use Proposition \ref{PropCompo}, one gets the result.\\
\noindent "$\Leftarrow$" Let $u$ be given by \eqref{u=volterra_C}, then 
\begin{align*}
u(t) = \sum_{k=0}^{n-1}\dfrac{\left(\varrho^k u\right)(a)}{k!}\big{(}\psi(t)-\psi(s)\big{)}^k  + \Iint{\apa,\psi}{a}f(t,u(t)).
\end{align*}
Apply the Riemann-Liouville derivative with respect to another function operator and use Proposition \ref{prop_semigpe_RL} yields
\begin{equation*}
\derg{\apa}{a}u(t) = \derg{\apa}{a}\left(\sum_{k=0}^{n-1}\dfrac{\left(\varrho^k u\right)(a)}{k!}\big{(}\psi(t)-\psi(s)\big{)}^k\right) + f(t,u(t)).
\end{equation*}

Then, the first equation of the system \eqref{fracds_c} follows from the definition of the fractional derivative with respect to another function in the Caputo sense. Finally, the initial conditions can be retrieved by applying the operator $\varrho^k = \catg{k}{a}$, $k=0,\, 1,\dots ,\, n-1$ (see Proposition \ref{prop_Cg=Igamma}) on both sides of \eqref{u=volterra_C} and using Lemma \ref{lem_der_C} and Proposition \ref{prop_semigpe_RL}.
\end{proof}

\begin{theo}\label{Theo_exist_unique}
Let $\apa > 0$ and $n = \lceil \apa \rceil$. Let $a_0, \dots, a_{n-1}$ in $\R$ and let $f : [a,b]\times \R \rightarrow \R$ be continuous function satisfying a Lipschitz condition with respect to its second variable. Then there exists a uniquely defined function $u\in AC^n_\psi[a,b]$ solving the initial value problem \eqref{fracds_c}.
\end{theo}
\begin{proof}
The function $F : (t,x) \mapsto f(\psi^{-1}(t),x)$ is continuous over $[\abar,\bbar]\times \R$ and fulfills a Lipschitz condition with respect to its second variable. It follows from  Theorem 6.8 in \cite{Diethelm10} that their exists a unique function $v \in AC^n[\abar,\bbar]$ solution of the system \eqref{fde1}. According to Theorem \ref{TheoEqCap} and Corollary \ref{CoroEq}, the function $u = v\circ \psi\in AC^n_\psi[a,b]$ is a solution of the system \eqref{fracds_c}. Finally, the uniqueness of $u$ is a direct consequence of the uniqueness of $v$.
\end{proof}

\begin{theo}\label{Theo_smoothness_1}
Let $\apa \in (0,1)$ and $\psi\in C^1[a,b]$ an increasing function such that $\psi'(t)\neq 0$ for all $t\in [a,b]$. Assume that
\begin{itemize}
\item[(A1)] The function $(t,x)\mapsto f(t,x)$ is of class $C^1$ over $[a,b]\times \R$.
\item[(A2)] The function $(t,x)\mapsto \partial_x f(t,x)$ is locally Lipschitz continuous in $x$.
\item[(A3)] The unique continuous solution of \eqref{fracds_c} exists on $[a,b]$.
\end{itemize}
Then the solution $u$ of \eqref{fracds_c} is of class $C[a,b]\cap C^1(a,b]$.
\end{theo}
\begin{proof}
The function $(t,x)\mapsto f(\psi^{-1}(t),x)$ is of class $C^1$ over $[\abar,\bbar]\times \R$ and the function $(t,x)\mapsto \partial_x f(\psi^{-1}(t),x)$ is locally Lipschitz continuous in $x$ for any $t\in [\abar,\bbar]$. From Theorem \ref{TheoEqCap} hereinabove and Theorem 1 in \cite{Miller71} we deduce that $v := u\circ \psi^{-1} \in C[\abar,\bbar]\cap C^1(\abar,\bbar]$, and hence $u\in C[a,b]\cap C^1(a,b]$.
\end{proof}

\section{Application: the Erd\'elyi-Kober operator }\label{SectionErd}
In this section, we derive several results related to the Erd\'elyi-Kober operators by using the concept of the fractional operators with respect to another function. First, we recall that the Erd\'elyi-Kober fractional integrals and derivatives (in the sense of Riemann-Liouville) are given respectively by
\begin{equation}\label{int_EK}
I^{\apa}_{a,\sigma,\eta}u(t) = \dfrac{\sigma\,t^{-\sigma(\apa+\eta)}}{\Gamma(\apa)}\int_a^t \dfrac{s^{\sigma\eta+\sigma-1}u(s)}{\left(t^\sigma-s^\sigma\right)^{1-\apa}}\,ds
\end{equation}
and
\begin{equation}\label{der_EK}
\Erd{\apa}{a,\sigma,\eta}u(t) = t^{-\sigma\eta} \left(\dfrac{1}{\sigma\,t^{\sigma-1}}\dfrac{d}{dt}\right)^n t^{\sigma(n+\eta)} I^{n-\apa}_{a,\sigma,\eta+\apa}u(t).
\end{equation}
If we set $\psi(t) = t^\sigma$, then one can check that the aforementioned operators can be written as
\begin{equation}\label{int_EK_g}
I^{\apa}_{a,\sigma,\eta}u(t) = t^{-\sigma(\apa+\eta)}\,I^{\apa,\psi}_{a}\left(t^{\sigma\eta}u(t)\right)
\end{equation}
and
\begin{equation}\label{der_EK_g}
\Erd{\apa}{a,\sigma,\eta}u(t) = t^{-\sigma\eta}\, \derg{\apa}{a}\left(t^{\sigma(\apa+\eta)}u(t)\right).
\end{equation}

Consequently, the Erd\'elyi-Kober fractional operators are closely related to the RL fractional operators with respect to another function. In the sequel, we shall consider $\psi(t)=t^\sigma$.
\begin{propos}
Let $\apa\in \R\backslash \N_0$ and $\eta \in \R\backslash \Z_0^-$. Let $u$ be a regular function. Then for any $t>0$
\begin{equation*}
\Erd{\apa}{0,\sigma,\eta}u(t) =\Gamma(1+\apa+\eta)\left(\dfrac{u(t)}{\Gamma(1+\eta)} + \sum_{k=1}^\infty \begin{pmatrix} \apa\\k\end{pmatrix}\dfrac{\sigma^{-k}}{\Gamma(1+k+\eta)}\sum_{j=1}^k \lma_{j,k}\,t^j\,u^{(j)}(t)\right),
\end{equation*}
where $\begin{pmatrix} \apa\\k\end{pmatrix} = \dfrac{\Gamma(\apa+1)}{k!\,\Gamma(\apa-k+1)}$ and $(\lma_{j,k})$ is the sequence given by
\begin{align*}
\lma_{j,k} = \left\{
\begin{array}{ll}
0 & \ \text{ if }\ \ j=0 \ \text{ or }\ j > k \\
1 & \ \text{ if }\ \ j = k \\
\lma_{j-1,k-1} + \big{(}j-(k-1)\sigma\big{)}\,\lma_{j,k-1} & \ \text{ if }\ \ 1\leq j < k.
\end{array}
\right.
\end{align*}
\end{propos}
\begin{proof}
The result can be proved using equation \eqref{der_EK_g}, Proposition \ref{prop_Leibniz}, Lemma \ref{lem_der} and Appendix \ref{SectionAppB}.
\end{proof}
\begin{lemma}\label{lem_expo}
Let $a\geq 0$, $\apa > 0$, $\sigma>0$ and $\eta$ and $\lma$ in $\R$. Let $p\geq 1$ and $u\in L^p(a,b)$ with either $a>0$ or $a=0$ and $\eta > -1+\frac{1}{p\sigma}$. Then for a.e. $t\in [a,b]$ we have
\begin{equation*}
I^{\apa}_{a,\sigma,\eta}\left(t^{\sigma\lma}u(t)\right) = t^{\sigma\lma}\, I^{\apa}_{a,\sigma,\eta+\lma}u(t)
\end{equation*}
and
\begin{equation*}
\Erd{\apa}{a,\sigma,\eta}\left(t^{\sigma\lma}u(t)\right) = t^{\sigma\lma}\, \Erd{\apa}{a,\sigma,\eta+\lma}u(t).
\end{equation*}
\end{lemma}
\begin{proof}
A direct consequence of \eqref{int_EK_g} and \eqref{der_EK_g} respectively.
\end{proof}
\begin{lemma}
Let $a\geq 0$, $\sigma>0$ and $\eta\in \R$. Let $\apa > 0$ and $n = \lceil \apa \rceil$. Then
\begin{align*}
& \Iint{\apa}{a,\sigma,\eta}\left(t^{-\sigma\eta}(t^\sigma-a^\sigma)^\bta\right) = \dfrac{\Gamma(1+\bta)}{\Gamma(1+\bta+\apa)}t^{-\sigma(\apa+\eta)}(t^\sigma-a^\sigma)^{\bta+\apa}, & \bta > -1,\\
& \Erd{\apa}{a,\sigma,\eta}\left(t^{-\sigma(\apa+\eta)}(t^\sigma-a^\sigma)^\bta\right) = \dfrac{\Gamma(1+\bta)}{\Gamma(1+\bta-\apa)}t^{-\sigma\eta}\left(t^\sigma-a^\sigma\right)^{\bta-\apa}, & \bta > -1,\\
& \Erd{\apa}{a,\sigma,\eta}\left(t^{-\sigma(\apa+\eta)}(t^\sigma-a^\sigma)^{\apa-i}\right) = 0, & i=1,\dots n,
\end{align*}
\end{lemma}
\begin{proof}
From \eqref{int_EK_g} and \eqref{der_EK_g} we have
\begin{equation*}
I^{\apa}_{a,\sigma,\eta}\left(t^{-\sigma\eta}u(t)\right) = t^{-\sigma(\apa+\eta)}\,I^{\apa,\psi}_{a}u(t)
\end{equation*}
and
\begin{equation*}
\Erd{\apa}{a,\sigma,\eta}\left(t^{-\sigma(\apa+\eta)}u(t)\right) = t^{-\sigma\eta}\, \derg{\apa}{a}u(t).
\end{equation*}
Using Lemma \ref{lem_der}, we deduce for $\bta>-1$
\begin{align*}
\Iint{\apa}{a,\sigma,\eta}\left(t^{-\sigma\eta}(t^\sigma-a^\sigma)^\bta\right) & = t^{-\sigma(\apa+\eta)}\,I^{\apa,\psi}_{a}(t^\sigma-a^\sigma)^\bta \\
& = \dfrac{\Gamma(1+\bta)}{\Gamma(1+\bta+\apa)}t^{-\sigma(\apa+\eta)}(t^\sigma-a^\sigma)^{\bta+\apa},
\end{align*}
and
\begin{align*}
\Erd{\apa}{a,\sigma,\eta}\left(t^{-\sigma(\apa+\eta)}(t^\sigma-a^\sigma)^\bta\right)& = t^{-\sigma\eta}\, \derg{\apa}{a}(t^\sigma-a^\sigma)^\bta\\
& = \dfrac{\Gamma(1+\bta)}{\Gamma(1+\bta-\apa)}t^{-\sigma\eta}\,(t^\sigma-a^\sigma)^{\bta-\apa}.
\end{align*}
The third assertion can be obtained similarly.
\end{proof}
\begin{lemma}
Let $a>0$, $\sigma>0$ and $\eta\in \R$. Let $\apa > 0$ and $n = \lceil \apa \rceil$. Then for $\bta/\sigma > -1$ and $t>a$
\begin{align*}
I^{\apa}_{a,\sigma,\eta} t^\bta & = t^{-\sigma(\apa+\eta)}\sum_{k\geq 0} \dfrac{\Gamma(1+\frac{\bta}{\sigma}+\eta)\,a^{\bta+\sigma(\eta-k)}}{\Gamma(\frac{\bta}{\sigma}+\eta-k)\,\Gamma(1+k+\apa)}(t^\sigma-a^\sigma)^{k+\apa} ,\\
\Erd{\apa}{a,\sigma,\eta} t^\bta & = t^{-\sigma\eta}\sum_{k\geq 0} \dfrac{\Gamma(1+\frac{\bta}{\sigma}+\apa+\eta)\ a^{\bta+\sigma(\apa+\eta-k)}}{\Gamma(\frac{\bta}{\sigma}+\apa+\eta-k)\,\Gamma(1+k-\apa)}(t^\sigma-a^\sigma)^{k-\apa}.
\end{align*}
In particular, for $m\in \N_0$ we have
\begin{align*}
I^{\apa}_{a,\sigma,\eta} t^m & = t^{-\sigma(\apa+\eta)}\sum_{k=0}^m \dfrac{\Gamma(1+\frac{m}{\sigma}+\eta)\,a^{m+\sigma(\eta-k)}}{\Gamma(\frac{m}{\sigma}+\eta-k)\,\Gamma(1+k+\apa)}(t^\sigma-a^\sigma)^{k+\apa} ,\\
\Erd{\apa}{a,\sigma,\eta} t^m & = t^{-\sigma\eta}\sum_{k=0}^m \dfrac{\Gamma(1+\frac{m}{\sigma}+\apa+\eta)\ a^{m+\sigma(\apa+\eta-k)}}{\Gamma(\frac{m}{\sigma}+\apa+\eta-k)\,\Gamma(1+k-\apa)}(t^\sigma-a^\sigma)^{k-\apa}.
\end{align*}
\end{lemma}
\begin{proof}
The results can be derived by taking $u(t) = t^\bta$ in \eqref{int_EK_g} and \eqref{der_EK_g} respectively, then using Lemma \ref{lem_der_power}.
\end{proof}

\begin{propos}\label{prop_permut_EK}
Let $a\geq 0$, $\apa > 0$, $\bta>0$ and $\sigma>0$. Let $p\geq 1$ and $u\in L^p(a,b)$ and let $\eta$ and $\lma$ be arbitrary real numbers if $a>0$ or $\min(\eta,\lma) > -1+\frac{1}{p\sigma}$ if $a=0$. Then, for a.e. $t\in [a,b]$
\begin{equation*}
\Iint{\apa}{a,\sigma,\eta}\Iint{\bta}{a,\sigma,\lma}u(t) = \Iint{\bta}{a,\sigma,\lma}\Iint{\apa}{a,\sigma,\eta}u(t).
\end{equation*}
\end{propos}
\begin{proof}
Using \eqref{int_EK_g} and Proposition \ref{prop_permut_RL} we obtain
\begin{align*}
\Iint{\apa}{a,\sigma,\eta}\Iint{\bta}{a,\sigma,\lma}u(t) & = t^{-\sigma(\apa+\eta)}\,I^{\apa,\psi}_{a}\left(t^{\sigma(\eta-\bta-\lma)}\,I^{\bta,\psi}_{a}\left(t^{\sigma\lma} u(t)\right)\right)\\
& = t^{-\sigma(\bta+\lma)}\,I^{\bta,\psi}_{a}\left(t^{-\sigma(\apa+\eta-\lma)}\,I^{\apa,\psi}_{a}\left(t^{\sigma\eta} u(t)\right)\right)\\
& = t^{-\sigma(\bta+\lma)}\,I^{\bta,\psi}_{a}\left(t^{\sigma\lma}\,I^{\apa}_{a,\sigma,\eta} u(t)\right)\\
& = \Iint{\bta}{a,\sigma,\lma}\Iint{\apa}{a,\sigma,\eta}u(t).
\end{align*}
\end{proof}
\begin{lemma}\label{II_EK_g}
Let $\apa>0$, $\bta>0$, $\sigma>0$, $p\geq 1$ and $u\in L^p(a,b)$ with either $a>0$ or $a=0$ and $\eta > -1+\frac{1}{p\sigma}$. Then for a.e. $t\in [a,b]$
\begin{equation*}
I^{\bta}_{a,\sigma,\eta}I^{\apa}_{a,\sigma,\eta+\bta}u(t) = I^{\apa}_{a,\sigma,\eta+\bta}I^{\bta}_{a,\sigma,\eta}u(t) = I^{\apa+\bta}_{a,\sigma,\eta}u(t) .
\end{equation*}
\end{lemma}
\begin{proof}
The first equality is a consequence of Proposition \ref{prop_permut_EK}. To obtain the second equality, we use \eqref{int_EK_g} and Proposition \ref{prop_semigpe_RL}. We have
\begin{align*}
I^{\apa}_{a,\sigma,\eta+\bta}I^{\bta}_{a,\sigma,\eta}u(t) & = I^{\apa}_{a,\sigma,\eta+\bta}\left(t^{-\sigma(\bta+\eta)}\,I^{\bta,\psi}_{a}\left(t^{\sigma\eta}u(t)\right)\right) \\
& = t^{-\sigma(\apa+\bta+\eta)}I^{\apa,\psi}_{a}\left(I^{\bta,\psi}_{a}\left(t^{\sigma\eta}u(t)\right)\right) \\
& = t^{-\sigma(\apa+\bta+\eta)}I^{\apa+\bta,\psi}_{a}\left(t^{\sigma\eta}u(t)\right)\\
& = I^{\apa+\bta}_{a,\sigma,\eta}u(t).
\end{align*}
\end{proof}
\begin{propos}
Assume the hypothesis of Lemma \ref{II_EK_g}. Then for $\apa \geq \bta$ and $u\in L^p(a,b)$
\begin{equation*}
\Erd{\bta}{a,\sigma,\eta}I^{\apa}_{a,\sigma,\eta+\bta-\apa}u(t) = I^{\apa-\bta}_{a,\sigma,\eta+\bta-\apa}u(t).
\end{equation*}
In particular
\begin{equation*}
\Erd{\apa}{a,\sigma,\eta}I^{\apa}_{a,\sigma,\eta}u(t) = u(t).
\end{equation*}
\end{propos}
\begin{proof}
Using \eqref{int_EK_g}, \eqref{der_EK_g} and Proposition \ref{prop_semigpe_RL}, we get
\begin{align*}
\Erd{\bta}{a,\sigma,\eta}I^{\apa}_{a,\sigma,\eta+\bta-\apa}u(t) & = t^{-\sigma\eta}\,\derg{\bta}{a}\left(I^{\apa,\psi}_{a}\left(t^{\sigma(\eta+\bta-\apa)}u(t)\right)\right) \\
& = t^{-\sigma\eta}\,I^{\apa-\bta,\psi}_{a}\left(t^{\sigma(\eta+\bta-\apa)}u(t)\right) \\
& = I^{\apa-\bta}_{a,\sigma,\eta+\bta-\apa}u(t).
\end{align*}
\end{proof}
\begin{propos}
Assume the hypothesis of Lemma \ref{II_EK_g}. Assume in addition that $t^{\sigma(n+\eta)}I^{n-\apa}_{a,\sigma,\apa+\eta}u\in AC^n_\psi[a,b]$, then
\begin{equation*}
I^{\apa}_{a,\sigma,\eta}\Erd{\apa}{a,\sigma,\eta}u(t) = u(t) - \sum_{k=1}^n \dfrac{\Erd{\apa-k}{a,\sigma,\eta+\apa}\left(t^{\sigma(k-\apa)}\,u(t)\right)(a)}{\Gamma(\apa-k+1)} \left(t^\sigma-a^\sigma\right)^{\apa-k}.
\end{equation*}
\end{propos}
\begin{proof}
First, let us notice that \eqref{int_EK_g} yields $I^{n-\apa,\psi}_{a}\left(t^{\sigma(\apa+\eta)}\,u\right) = t^{\sigma(n+\eta)}I^{n-\apa}_{a,\sigma,\apa+\eta}u \in AC^n_\psi[a,b]$. It follows using \eqref{int_EK_g}, \eqref{der_EK_g}, Proposition \ref{prop_IagDag} and Lemma \ref{lem_expo}
\begin{align*}
\hspace*{-5mm} I^{\apa}_{a,\sigma,\eta}\Erd{\apa}{a,\sigma,\eta}u(t) & = t^{-\sigma(\apa+\eta)}\,I^{\apa,\psi}_{a}\left(t^{\sigma\eta}\,\Erd{\apa}{a,\sigma,\eta}u(t)\right)\\
& = t^{-\sigma(\apa+\eta)}\,I^{\apa,\psi}_{a}\,\derg{\apa}{a}\left(t^{\sigma(\apa+\eta)}u(t)\right)\\
& = t^{-\sigma(\apa+\eta)}\left[t^{\sigma(\apa+\eta)}u(t) - \sum_{k=1}^n \dfrac{\derg{\apa-k}{a}\left(t^{\sigma(\apa+\eta)}\,u(t)\right)(a)}{\Gamma(\apa-k+1)} \left(t^\sigma-a^\sigma\right)^{\apa-k}\right]\\
& = u(t) - \sum_{k=1}^n \dfrac{\Erd{\apa-k}{a,\sigma,\eta+\apa}\left(t^{\sigma(k-\apa)}\,u(t)\right)(a)}{\Gamma(\apa-k+1)} \left(t^\sigma-a^\sigma\right)^{\apa-k}.
%
\end{align*}
\end{proof}

\begin{propos}
Let $\apa>0$ and $\bta>0$ such that $n-1< \apa \leq n$ and $m-1< \bta \leq m$ with $n, \, m\in \N$. Let $u\in L^p(a,b)$, $p\geq 1$ with either $a>0$ or $a=0$ and $\eta > -1+\frac{1}{p\sigma}$. Assume in addition $t^{\sigma(m+\apa+\eta)}\Iint{m-\bta}{a,\sigma,\apa+\bta+\eta} u \in AC^m_\psi[a,b]$, then for a.e. $t\in [a,b]$
\begin{align*}
\Erd{\apa}{a,\sigma,\eta}\left(\Erd{\bta}{a,\sigma,\eta+\apa}u\right)(t) & = \Erd{\apa+\bta}{a,\sigma,\eta}u(t) \\
& \quad - \sum_{k=1}^{m-1} \dfrac{\Erd{\bta-k}{a,\sigma,\eta}\left(t^{\sigma(\apa+k)}u(t)\right)(a)}{\Gamma(1-k-\apa)} \, (t^\sigma-a^\sigma)^{-k-\apa}\\
& \quad - \dfrac{(t^\sigma-a^\sigma)^{-m-\apa}}{\Gamma(1-m-\apa)} \lim_{t\to a}\Iint{m-\bta}{a,\sigma,\eta-m+\bta}\left(t^{\sigma(\apa+m)}u(t)\right).
\end{align*}
\end{propos}
\begin{proof}
Using \eqref{int_EK_g}, \eqref{der_EK_g} and Proposition \ref{prop_Da_Db}, one obtain
\begin{align*}
& \Erd{\apa}{a,\sigma,\eta}\left(\Erd{\bta}{a,\sigma,\eta+\apa}u\right)(t) \\
& = t^{-\sigma\eta}\,\derg{\apa}{a}\left(t^{\sigma(\apa+\eta)}\,\Erd{\bta}{a,\sigma,\eta+\apa}u(t)\right) \\
& = t^{-\sigma\eta}\,\derg{\apa}{a}\,\derg{\bta}{a}\left(t^{\sigma(\apa+\bta+\eta)}u(t)\right) \\
& = t^{-\sigma\eta}\left\{\derg{\apa+\bta}{a}\left(t^{\sigma(\apa+\bta+\eta)}u(t)\right)  - \sum_{k=1}^{m-1} \dfrac{\derg{\bta-k}{a}\left(t^{\sigma(\apa+\bta+\eta)}u(t)\right)(a)}{\Gamma(1-k-\apa)}(t^\sigma-a^\sigma)^{-k-\apa} \right.\\
& \qquad \qquad \left. - \dfrac{(t^\sigma-a^\sigma)^{-m-\apa}}{\Gamma(1-m-\apa)}\,\lim_{t\to a}\Iint{m-\bta,\psi}{a}\left(t^{\sigma(\apa+\bta+\eta)}u(t)\right)\right\}\\
& = \Erd{\apa+\bta}{a,\sigma,\eta}u(t) - \sum_{k=1}^{m-1} \dfrac{\Erd{\bta-k}{a,\sigma,\eta}\left(t^{\sigma(\apa+k)}u(t)\right)(a)}{\Gamma(1-k-\apa)}\,(t^\sigma-a^\sigma)^{-k-\apa}\\
& \quad - \dfrac{(t^\sigma-a^\sigma)^{-m-\apa}}{\Gamma(1-m-\apa)}\,\lim_{t\to a}\Iint{m-\bta}{a,\sigma,\eta-m+\bta}\left(t^{\sigma(\apa+m)}u(t)\right) .
\end{align*}
\end{proof}
\begin{propos}\label{prop_I-beta_EK}
Let $a>0$, $\apa>0$ with $\apa\not\in \N$ and $n=\lceil \apa \rceil$. Then for any $u\in C^n[a,b]$ and any $t\in [a,b]$
\begin{equation*}
\Erd{\apa}{a,\sigma,\eta} u(t) = \Iint{-\apa}{a,\sigma,\eta+\apa} u(t).
\end{equation*}
\end{propos}
\begin{proof}
Since $a>0$, then $\psi(t) = t^\sigma\in C^n[a,b]$ and $t^{\sigma(\apa+\eta)}u(t)\in C^n[a,b]$. It follows from \eqref{int_EK_g}, \eqref{der_EK_g} and Proposition \ref{prop_I-beta} that
\begin{align*}
\Iint{-\apa}{a,\sigma,\eta+\apa} u(t) & = t^{-\sigma\eta}\Iint{-\apa,\psi}{a}\left(t^{\sigma(\apa+\eta)} u(t)\right) \\
& = t^{-\sigma\eta}\,\derg{\apa}{a}\left(t^{\sigma(\apa+\eta)} u(t)\right)\\
& = \Erd{\apa}{a,\sigma,\eta}u(t).
\end{align*}
\end{proof}
\begin{rem}
In case $a=0$, then Proposition \ref{prop_I-beta_EK} holds also true if one suppose in addition that $\sigma\geq n$ and $t^{\sigma(\apa+\eta)}u\in C^n[0,b]$.
\end{rem}
Let $a\geq 0$, $\apa>0$ and $n= \lceil \apa \rceil$. Consider the following fractional differential system
\begin{align}\label{fracds_EK}
\left\{
\begin{array}{l}
\Erd{\apa}{a,\sigma,\eta}u(t) = f(t,u(t)), \quad t\in [a,b], \\ \\
\Erd{\apa-k}{a,\sigma,0}\left(t^{\sigma(k+\eta)}u(t)\right)(a) = a_k, \ \ k=1,\dots ,\, n-1, \ \ \spl{\lim_{t\to a}}\Iint{n-\apa}{a,\sigma,\apa-n}\left(t^{\sigma(n+\eta)}u(t)\right) = a_n.
\end{array}
\right.
\end{align}
with $f$ a given function and $a_k\in \R$ for $k= 1,\dots ,\, n$. Then, we have the following result.
\begin{theo}\label{Theo_volt_EK}
Let $f : (a,b]\times \R \rightarrow \R$ be a function such that the mapping $t\mapsto f(t,t^{-\sigma(\apa+\eta)}x)\in L^1_{\psi^{\eta+1}}(a,b)$ for any $x\in \R$\textnormal{(}\footnote{We recall $L^1_{\psi^\lma}(a,b) = \left\{f:[a,b]\rightarrow \Cx,\ \int_a^b s^{\sigma\lma-1}|f(s)|\,ds < \infty\right\}$ (see equation \eqref{Lpg}).}\textnormal{)}. Then a function $u\in L^1_{\psi^{\apa+\eta+1}}(a,b)$ is a solution of \eqref{fracds_EK} if and only if $u$ is a solution of the nonlinear Volterra integral equation of the second kind
\begin{align}\label{u=volterra_EK}
u(t) = t^{-\sigma(\apa+\eta)}\,\sum_{k=1}^{n}\dfrac{a_k}{\Gamma(\apa-k+1)}\big{(}t^\sigma-a^\sigma\big{)}^{\apa-k}  + \Iint{\apa}{a,\sigma,\eta}f(t,u(t)),
\end{align}
with $a_k\!=\!\Erd{\apa-k}{a,\sigma,0}\left(t^{\sigma(k+\eta)}u(t)\right)(a)$, $k=1,\dots ,\, n-1$, and $a_n\!=\!\spl{\lim_{t\to a}}\Iint{n-\apa}{a,\sigma,\apa-n}\left(t^{\sigma(n+\eta)}u(t)\right)$.
\end{theo}
\begin{proof}
First, notice that we have by equation \eqref{int_EK_g} and \eqref{der_EK_g} 
$$\Erd{\apa-k}{a,\sigma,0}\left(t^{\sigma(k+\eta)}u(t)\right) = \derg{\apa-k}{a}\left(t^{\sigma(\apa+\eta)}u(t)\right)$$ and
$$\Iint{n-\apa}{a,\sigma,\apa-n}\left(t^{\sigma(n+\eta)}u(t)\right) = \Iint{n-\apa,\psi}{a}\left(t^{\sigma(\apa+\eta)}u(t)\right).$$
Using again equations \eqref{int_EK_g} and \eqref{der_EK_g}, one may deduce that $u$ is a solution of \eqref{fracds_EK} if and only if $v(t) = t^{\sigma(\apa+\eta)}u(t)$ is a solution of the system
\begin{align}\label{fracds_EK_equiv}
\left\{
\begin{array}{l}
\derg{\apa}{a}v(t) = F(t,v(t)), \quad t\in [a,b] \\ \\
\derg{\apa-k}{a}v(a) = a_k, \ \ k=1,\dots ,\, n-1, \ \ \spl{\lim_{t\to a}}\Iint{n-\apa,\psi}{a}v(t) = a_n,
\end{array}
\right.
\end{align}
with $F(t,x) = t^{\sigma\eta}\,f(t,t^{-\sigma(\apa+\eta)}x)$. From the hypothesis, $t\mapsto F(t,\cdot)\in L^1_\psi(a,b)$. It follows using Proposition \ref{prop_volt}
\begin{equation*}
v(t) = \sum_{k=1}^{n}\dfrac{a_k}{\Gamma(\apa-k+1)}\big{(}t^\sigma-a^\sigma\big{)}^{\apa-k}  + \dfrac{\sigma}{\Gamma(\apa)}\int_a^t \dfrac{s^{\sigma-1}F(s,v(s))}{(t^\sigma-s^\sigma)^{1-\apa}}\,ds.
\end{equation*}
Multiplying $v(t)$ by $t^{-\sigma(\apa+\eta)}$ and substituting $F$ by its expression yields \eqref{u=volterra_EK}. Finally, since $v\in L^1_\psi(a,b)$ then $u\in L^1_{\psi^{\apa+\eta+1}}(a,b)$, and the proof is completed.
\end{proof}

Define the space
$$L^\apa_{\sigma,\eta}(a,b) := \left\{\varphi\in L^1_{\psi^{\apa+\eta+1}}(a,b),\ \Erd{\apa}{a,\sigma,\eta}\varphi\in L^1_{\psi^{\eta+1}}(a,b)\right\},$$
with $\psi(t)=t^\sigma$. Then, we have the following result.
\begin{theo}\label{Theo_carathe}
Let $\apa>0$ and $n=\lceil\apa\rceil$. Let $f : (a,b]\times \R \rightarrow \R$ be a function such that $t\mapsto f(t,t^{-\sigma(\apa+\eta)}x)\in L^1_{\psi^{\eta+1}}(a,b)$ for any $x\in \R$. Assume their exists $K>0$ such that for all $t\in (a,b]$ and for all $(x,y)\in \R^2$
\begin{equation}\label{carathe}
\left|f(t,x)-f(t,y)\right| \leq K\,t^{\sigma\apa}\left|x-y\right|.
\end{equation}
Then the Cauchy problem \eqref{fracds_EK} admits a unique solution $u\in L^\apa_{\sigma,\eta}(a,b)$.
\end{theo}
\begin{proof}
We have established in Theorem \ref{Theo_volt_EK} that $u(t)$ is a solution of \eqref{fracds_EK} if and only if $v(t):=t^{\sigma(\apa+\eta)}u(t)$ is a solution of \eqref{fracds_EK_equiv}. On another hand, a straightforward computation shows that if $f$ satisfies \eqref{carathe} then $(t,x)\mapsto F(t,x):=t^{\sigma\eta}\,f(t,t^{-\sigma(\apa+\eta)}x)$ is Lipschitzian with respect to its second variable. Using Theorem \ref{exist_uniq_RL}, we deduce the existence and the uniqueness of the solution $v\in L^\apa_\psi(a,b)$ of \eqref{fracds_EK_equiv}, where $L^\apa_\psi(a,b)$ is given by \eqref{Lalphag}. It follows that $u(t) = t^{-\sigma(\apa+\eta)}v(t)\in L^1_{\psi^{\apa+\eta+1}}(a,b)$ and by using \eqref{der_EK_g} that $ \Erd{\apa}{a,\sigma,\eta}u\in L^1_{\psi^{\eta+1}}(a,b)$. Thus, $u\in L^\apa_{\sigma,\eta}(a,b)$ is the unique solution of \eqref{fracds_EK}.
\end{proof}
\begin{theo}\label{Theo_EK_equiv}
Assume the hypothesis of Theorem \ref{Theo_carathe}. Then $u$ is a solution of \eqref{fracds_EK} if and only if $u$ can be written as $u(t)=t^{-\sigma(\apa+\eta)}\omega(t^\sigma)$, where $\omega \in L^1(a^\sigma,b^\sigma)$ is the unique solution of the system
\begin{align*}
\left\{
\begin{array}{l}
\der{\apa}{a^\sigma}\omega(t) = G(t,\omega(t)), \quad t\in [a^\sigma,b^\sigma] \\ \\
\der{\apa-k}{a^\sigma}\omega(a^\sigma) = a_k, \quad k=1,\dots ,\, n-1, \ \ \spl{\lim_{t\to a^\sigma}}\Iint{n-\apa}{a^\sigma}\omega(t) = a_n.
\end{array}
\right.
\end{align*}
with $G(t,x) = t^{\eta} f(t^{1/\sigma},t^{-(\apa+\eta)}x)$.
\end{theo}
\begin{proof}
According to Theorem \ref{Theo_volt_EK}, $u$ is a solution of \eqref{fracds_EK} if and only if $v(t):=t^{\sigma(\apa+\eta)}u(t)$ is a solution of \eqref{fracds_EK_equiv} with $F(t,x):=t^{\sigma\eta}\,f(t,t^{-\sigma(\apa+\eta)}x)$. Now, using Corollary \ref{CoroEq_RL} we have that $v$ is a solution of \eqref{fracds_EK_equiv} if and only if $v=\omega\circ \psi$ where $\omega\in L^1(a^\sigma,b^\sigma)$ is the solution of the system
\begin{align*}
\left\{
\begin{array}{l}
\der{\apa}{a^\sigma}\omega(t) = G(t,\omega(t)), \quad t\in [a^\sigma,b^\sigma] \\ \\
\der{\apa-k}{a^\sigma}\omega(a^\sigma) = a_k, \quad k=1,\dots ,\, n-1, \ \ \spl{\lim_{t\to a^\sigma}}\Iint{n-\apa}{a^\sigma}\omega(t) = a_n.
\end{array}
\right.
\end{align*}
with $G(t,x) = F(t^{1/\sigma},x)=t^{\eta} f(t^{1/\sigma},t^{-(\apa+\eta)}x)$. This achieves the proof.
\end{proof}
\begin{rem}
For more results in concern with the Erd\'elyi-Kober operators in case $a=0$, we refer the reader to the book of Kiryakova \cite{Kiryakova94} and the references therein.
\end{rem}

Now, we give some explicit examples to illustrate our ideas.\\

\textbf{Example 1:} Let $\psi\in C^1(a,+\infty)$ with $\psi'>0$ and consider the following Volterra integral equation of first kind
\begin{equation}\label{Volt_first}
\int_a^t \dfrac{e^{\lma(t-s)}u(s)}{\sqrt{\psi(t)-\psi(s)}}\,ds = f(t).
\end{equation}
Equation \eqref{Volt_first} can be written as 
\begin{equation}\label{Volt_first_modif}
\Gamma\left(1/2\right)\Iint{1/2,\psi}{a}\left(\dfrac{e^{-\lma t}u(t)}{\psi'(t)}\right) = e^{-\lma t}f(t).
\end{equation}
Apply $\derg{1/2}{a}$ on both sides of \eqref{Volt_first_modif} and use Proposition \ref{prop_semigpe_RL} yields $$u(t) = \dfrac{1}{\Gamma\left(1/2\right)} e^{\lma t} \psi'(t)\derg{1/2}{a}\left(e^{-\lma t}f(t)\right),$$
or equivalently
\begin{equation}\label{sol_Volt_first}
u(t) = \dfrac{1}{\pi} \dfrac{d}{dt}\int_a^t \dfrac{e^{\lma(t-s)}\psi'(s)f(s)}{\sqrt{\psi(t)-\psi(s)}}\,ds.
\end{equation}
This result is in accordance with a similar solution given in \cite{Pol08}.\\

\textbf{Example 2:} Let $\apa\in(0,1)$, $\lma\in \R$, $\psi\in C^1(a,+\infty)$ with $\psi'>0$, and $v\in C(a,+\infty)$ with $v(t)\neq 0$ for all $t$. Let consider the following Volterra integral equation of second kind
\begin{equation}\label{Volt_sec}
u(t) + \lma\dfrac{\psi'(t)}{v(t)}\int_a^t \dfrac{v(s)u(s)}{(\psi(t)-\psi(s))^{1-\apa}}\,ds = f(t),
\end{equation}
which can be written as 
\begin{equation*}
u(t) + \lma\Gamma\left(\apa\right)\dfrac{\psi'(t)}{v(t)}\Iint{\apa,\psi}{a}\left(\dfrac{vu}{\psi'}\right)(t) = f(t),
\end{equation*}
or equivalently by using Proposition \ref{prop_IDg}
\begin{equation*}
u(t) + \lma\Gamma\left(\apa\right)\dfrac{\psi'(t)}{v(t)}\Iint{\apa}{\abar}\left(\dfrac{vu}{\psi'}\circ\psi^{-1}\right)(\psi(t)) = f(t).
\end{equation*}
Denote $x=\psi(t)$ and $y = \dfrac{vu}{\psi'}\circ\psi^{-1}$, we obtain after simplification
\begin{equation*}
y(x) + \lma\Gamma\left(\apa\right)\Iint{\apa}{\abar}y(x) = F(x) 
\end{equation*}
with $F = \dfrac{f v}{\psi'}\circ \psi^{-1}$. Using the classical theory for Volterra integral equations, one gets
\begin{equation*}
y(x) = G(x) + \Lambda\int_{\abar}^x e^{\Lambda(x-\tau)}G(\tau)d\tau,
\end{equation*}
with $\Lambda = \lma^2\Gamma\left(\apa\right)\Gamma\left(1-\apa\right)$ and $G(x) = F(x) - \lma\Gamma\left(\apa\right)I^{\apa}_{\abar}F(x)$. Finally, the solution of \eqref{Volt_sec} is given by
\begin{equation}\label{sol_Volt_sec}
u(t) = \dfrac{\psi'(t)}{v(t)} \left(\Xi(t) + \Lambda\int_a^t e^{\Lambda(\psi(t)-\psi(s))}\,\Xi(s)\,\psi'(s)\,ds\right)
\end{equation}
with $\Xi(t)=\frac{v(t)f(t)}{\psi'(t)} - \lma \int_a^t\frac{v(s)f(s)ds}{(\psi(t)-\psi(s))^{1-\apa}}$. For instance, if $w$ is a given function, then the solution of the following integral equation
\begin{equation*} u(t) + \dfrac{\lma}{w(t)} \int_a^t \dfrac{w(s)\psi'(s)u(s)}{\sqrt{\psi(t)-\psi(s)}}\,ds = f(t)
\end{equation*}
can be found by considering $v(t) = w(t)\psi'(t)$ in \eqref{Volt_sec}, yielding
\begin{equation*}
u(t) = \dfrac{1}{w(t)} \left(\Xi(t) + \pi\lma^2\int_a^t e^{\pi\lma^2(\psi(t)-\psi(s))}\,\Xi(s)\,\psi'(s)\,ds\right),
\end{equation*}
with $\Xi(t)=w(t)f(t)-\lma \spl{\int_a^t\frac{w(s)\psi'(s)f(s)ds}{\sqrt{\psi(t)-\psi(s)}}}$.\\

\textbf{Example 3:} Let $\apa\in(0,1)$, $\sigma>0$, $\eta\in\R$, $\lma\in \R$ and $u_0\in \R$. Consider the system
\begin{align}\label{Ex_EK}
\left\{
\begin{array}{l}
\Erd{\apa}{0,\sigma,\eta}u(t) = \lma\,t^{\sigma\apa}u(t),  \ \ t > 0,\\ \\
\lim_{t\to 0}\Iint{1-\apa}{0,\sigma,\apa-1}\left(t^{\sigma(1+\eta)}u(t)\right) = u_0.
\end{array}
\right.
\end{align}
In view of Theorem \ref{Theo_EK_equiv}, the solution of \eqref{Ex_EK} writes $u(t)=t^{-\sigma(\apa+\eta)}\omega(t^\sigma)$, where $\omega$ is the solution of the system
\begin{align}\label{Ex_EK_eqv}
\left\{
\begin{array}{l}
\der{\apa}{0}\omega(t) = \lma\,\omega(t), \quad t>0 \\ \\
\lim_{t\to 0}\Iint{1-\apa}{0}\omega(t) = u_0.
\end{array}
\right.
\end{align}
The solution of \eqref{Ex_EK_eqv} is given by $\omega(t) = u_0\,e_\apa^{\lma t}$ (see e.g. \cite[eq (2.1.56)]{Kilbas06}) where 
\begin{equation}\label{e_alpha}
e_\apa^{\lma x} = x^{\apa-1}\sum_{k\geq 0}\dfrac{(\lma\,x^\apa)^k}{\Gamma((k+1)\apa)}
\end{equation} is the $\apa$-exponential function. It follows that the solution of \eqref{Ex_EK} is given by $$u(t) = u_0\,t^{-\sigma(\apa+\eta)}\,e_\apa^{\lma t^\sigma}.$$


\section{Numerical methods for the fractional derivatives with respect to another function}\label{SectionNum}
We aim at deriving high order numerical methods able to accurately approach the solutions of  systems involving the fractional derivative operator with respect to another function. Since the scaling function $\psi$, as well as the solution of the fractional system and/or its first derivative might not be smooth at the lower terminal $a$, the convergence rate of any numerical scheme could drastically be impaired (see e.g. Lemma \ref{lem_der}). Yet, some numerical methods have been introduced in the literature to deal with the fractional derivative operators with respect to another function \cite{Alm17,Bal17}. In \cite{Zen17}, a finite difference scheme of order $1-\apa$ is presented to solve linear systems involving the generalized fractional derivatives. Another scheme with higher convergence rate is given in \cite{Odi20}. On another hand, several accurate schemes with high convergence orders are available for the  Caputo operators \cite{Cai20,Diethelm97,Lin2007,Lubich83,Mokh20}. Using these latter in combination with the results established in the previous sections, we claim that one can adequately obtain optimal convergence orders schemes to numerically solve the fractional systems involving the integral or derivative operators with respect to another function.
\subsection{Example 1}\label{exmp1}
\begin{figure}
\begin{center}
\begin{tabular}{c}
\includegraphics[scale=0.3]{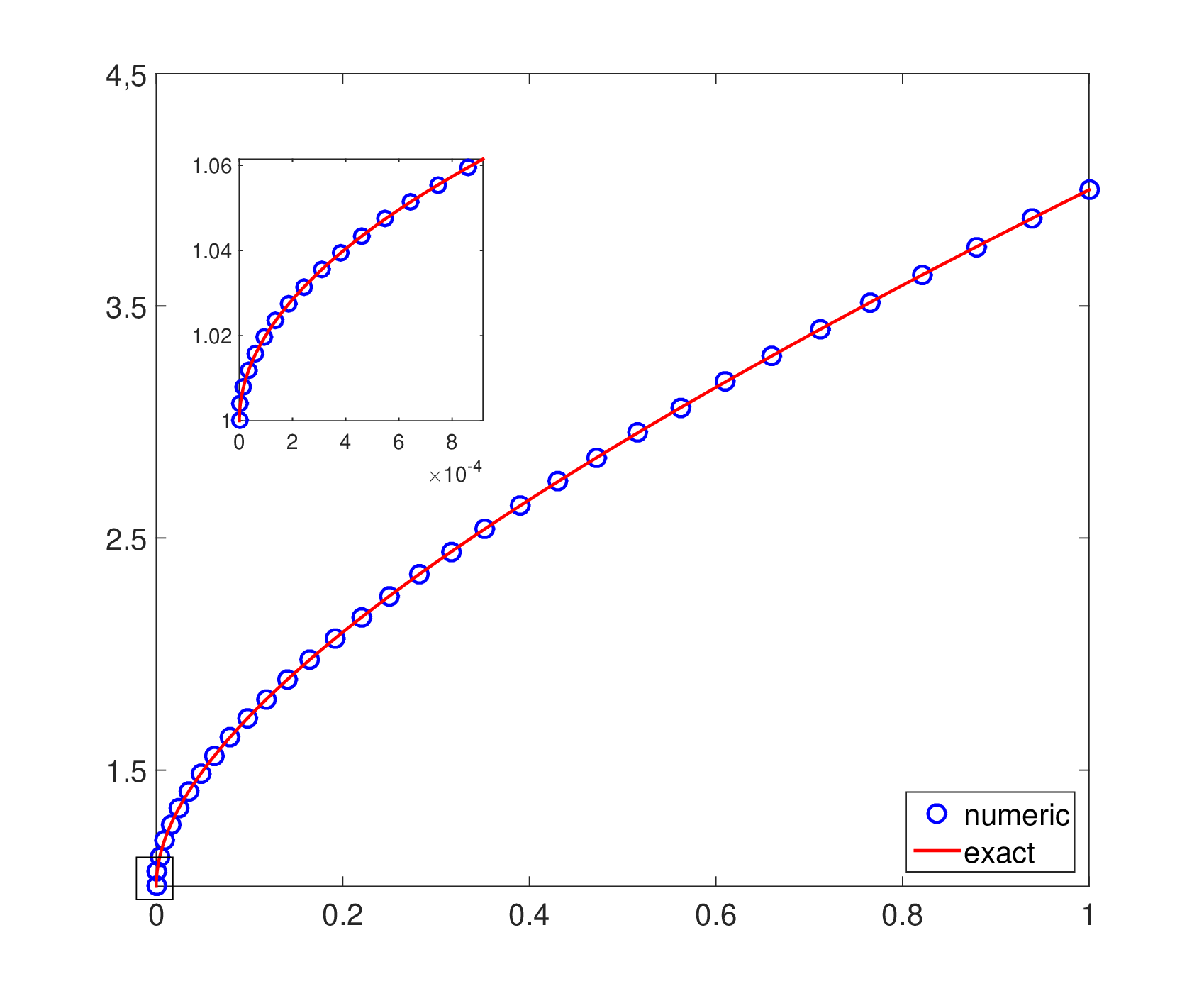}
\end{tabular}
\caption{Comparison between the numerical solution (blue circles) and the exact solution (red solid line) for the system \eqref{Example_1_sys} involving the Caputo derivative operator with respect to $\psi(t)=1+\sqrt{t}$.}\label{FigEx1}
\end{center}
\end{figure}

Let $\apa\in (0,1)$ and $\psi(t) = 1+\sqrt{t}$. Consider the system
\begin{align}\label{Example_1_sys}
\left\{
\begin{array}{ll}
\catg{\apa}{0}u(t) = \dfrac{2\,t^\frac{1-\apa}{2}}{\Gamma(3-\apa)}\left(2-\apa + \dfrac{\sqrt{t}}{(1+\sqrt{t})^2}\,u(t)\right), & \ t\geq 0,\\
u(0) = 1. &
\end{array}
\right.
\end{align}
In view of Theorem \ref{TheoEqCap}, the system \eqref{Example_1_sys} is equivalent to
\begin{align}\label{Example_1_syseqv}
\left\{
\begin{array}{l}
\cat{\apa}{1}v(t) = \dfrac{2\,(t-1)^{1-\apa}}{\Gamma(3-\apa)}\left(2-\apa + \dfrac{t-1}{t^2} \,v(t)\right), \quad t\geq 1, \\
v(1) = 1,
\end{array}
\right.
\end{align}
with $v = u\circ \psi^{-1}$. Using Lemma \ref{lem_der_C_power}, one may check that the solution of \eqref{Example_1_syseqv} is $v(t) = t^2$, yielding $u(t) = (1+\sqrt{t})^2$. Now, we aim at constructing a numerical method that approximates the solution $u$ of \eqref{Example_1_sys}. In view of the singularity of $u'$ at the origin, one can expect that the convergence rate of any (classical) numerical method would be impaired if directly applied to the system \eqref{Example_1_sys}. On the other hand, the solution $v$ of \eqref{Example_1_syseqv} is smooth and could be accurately approached using any suitable numerical method for the Caputo operator. An approximation of $u$ can then be obtained by a scaling technique. We choose (for instance) a finite difference method introduced in \cite{Lin2007}. This method is of order $2-\apa$ for sufficiently smooth solutions. Table \ref{Example_1_table} shows the $\ell^\infty$-errors for various values of $\apa$, and Figure \ref{FigEx1} displays a comparison between the solutions. We notice that the convergence orders are close to the theoretical value, in spite that $u$ is not sufficiently smooth at the origin. Hence, considering the equivalent systems instead of the original ones turns out to be advantageous especially in case of singular solutions approximation, since as well known, a specific care should be taken to address such lack of accuracy, by using for instance a suitable graded meshing \cite{Rice1969} or by plugging some non polynomial functions into the numerical solution in order to mimic the exact solution's singularity \cite{Cao2003,Ford2015}. Our method allows us to directly tackle such problems without the use of a specific mesh transform neither an additional computational cost provided the equivalent system's solutions are sufficiently regular.

\begin{table}
\begin{center}
\setlength\extrarowheight{4pt}
\begin{tabular}{@{\extracolsep{\fill}}|c|c|c||c|c||c|c|}
\cline{2-7}
\multicolumn{1}{c|}{}  & \multicolumn{2}{c||}{$\apa = 0.8$} & \multicolumn{2}{c||}{$\apa = 0.5$} & \multicolumn{2}{c|}{$\apa = 0.2$} \\ 
\cline{1-7} 
$N$ & $\|e\|_{\ell^\infty}$ & Order & $\|e\|_{\ell^\infty}$ & Order & $\|e\|_{\ell^\infty}$ & Order \\ 
\cline{1-7} 
$2^7$ & $7.8386\ (-04)$ & -- & $1.1891\ (-04)$ & -- & $1.2042\ (-05)$ & -- \\ 
\cline{1-7} 
$2^8$ & $3.4191\ (-04)$ & $1.1970$ & $4.2423\ (-05)$ & $1.4869$ & $3.5870\ (-06)$ & $1.7473$ \\ 
\cline{1-7} 
$2^9$ & $1.4901\ (-04)$ & $1.1982$ & $1.5094\ (-05)$ & $1.4909$ & $1.0622\ (-06)$ & $1.7557$ \\ 
\cline{1-7} 
$2^{10}$ & $6.4908\ (-05)$ & $1.1989$ & $5.3601\ (-06)$ & $1.4936$ & $3.1308\ (-07)$ & $1.7625$ \\ 
\cline{1-7} 
$2^{11}$ & $2.8265\ (-05)$ & $1.1994$ & $1.9009\ (-06)$ & $1.4956$ & $9.1916\ (-08)$ & $1.7681$ \\ 
\cline{1-7} 
$2^{12}$ & $1.2306\ (-05)$ & $1.1996$ & $6.7354\ (-07)$ & $1.4969$ & $2.6898\ (-08)$ & $1.7728$ \\ 
\hline
\end{tabular}
\vspace*{2mm}
\caption{Errors and convergence orders relative to system \eqref{Example_1_sys} for various values of $\apa$. Despite $u$ is not smooth at the origin, the convergence orders are close to the optimal theoretical value $2-\apa$, for all the choices of $\apa$.}\label{Example_1_table}
\end{center}
\end{table}

\subsection{Example 2}\label{exmp2}
Let $\apa\in (0,1)$ and $u_0$ and $\lma$ in $\R$. Consider the system
\begin{equation}\label{Example_2_sys}
\left\{
\begin{array}{ll}
\Had{\apa}{1}u(t) = \lma\left(\log t\right)^{1-\apa}u(t), & t\geq 1 \\
u(1) = u_0, &
\end{array}
\right.
\end{equation}
where $\Had{\apa}{1}$ stands for the Hadamard fractional derivative operator in the sense of Caputo, given by \cite{Kilbas06}
\begin{equation*}
\Had{\apa}{a}u(t) = \dfrac{1}{\Gamma(n-\apa)}\int_a^t \left(\log \frac{t}{\tau}\right)^{n-\apa-1}\left[\left(\tau\dfrac{d}{d\tau}\right)^nu\right](\tau)\,\dfrac{d\tau}{\tau}.
\end{equation*}
The system \eqref{Example_2_sys} can be written in term of the fractional derivative with respect to the function $\psi(t):=\log t$ as
\begin{equation*}
\left\{
\begin{array}{ll}
\catg{\apa}{1}u(t) = \lma\left(\log t\right)^{1-\apa}u(t), & t\geq 1 \\
u(1) = u_0, &
\end{array}
\right.
\end{equation*}
or equivalently by Theorem \ref{TheoEqCap}
\begin{equation}\label{Example_2_syseqv}
\left\{
\begin{array}{ll}
\cat{\apa}{0}v(t) = \lma\, t^{1-\apa}v(t), & t\geq 0 \\
v(0) = u_0, &
\end{array}
\right.
\end{equation}
with the relation $v = u\circ\psi^{-1}$. A straightforward computation shows that $v(t) = u_0\sum_{k\geq 0}a_k\left(\lma t\right)^k$ with $a_0=1$ and $a_k = \prod_{j=1}^k \frac{\Gamma(j+1-\apa)}{\Gamma(j+1)}$, yielding
\begin{equation*}
u(t) = u_0 + \sum_{k\geq 1}\prod_{j=1}^k \frac{\Gamma(j+1-\apa)}{\Gamma(j+1)}\left(\lma \log t\right)^k.
\end{equation*}

\begin{figure}[h!]
\begin{center}
\includegraphics[scale=0.3]{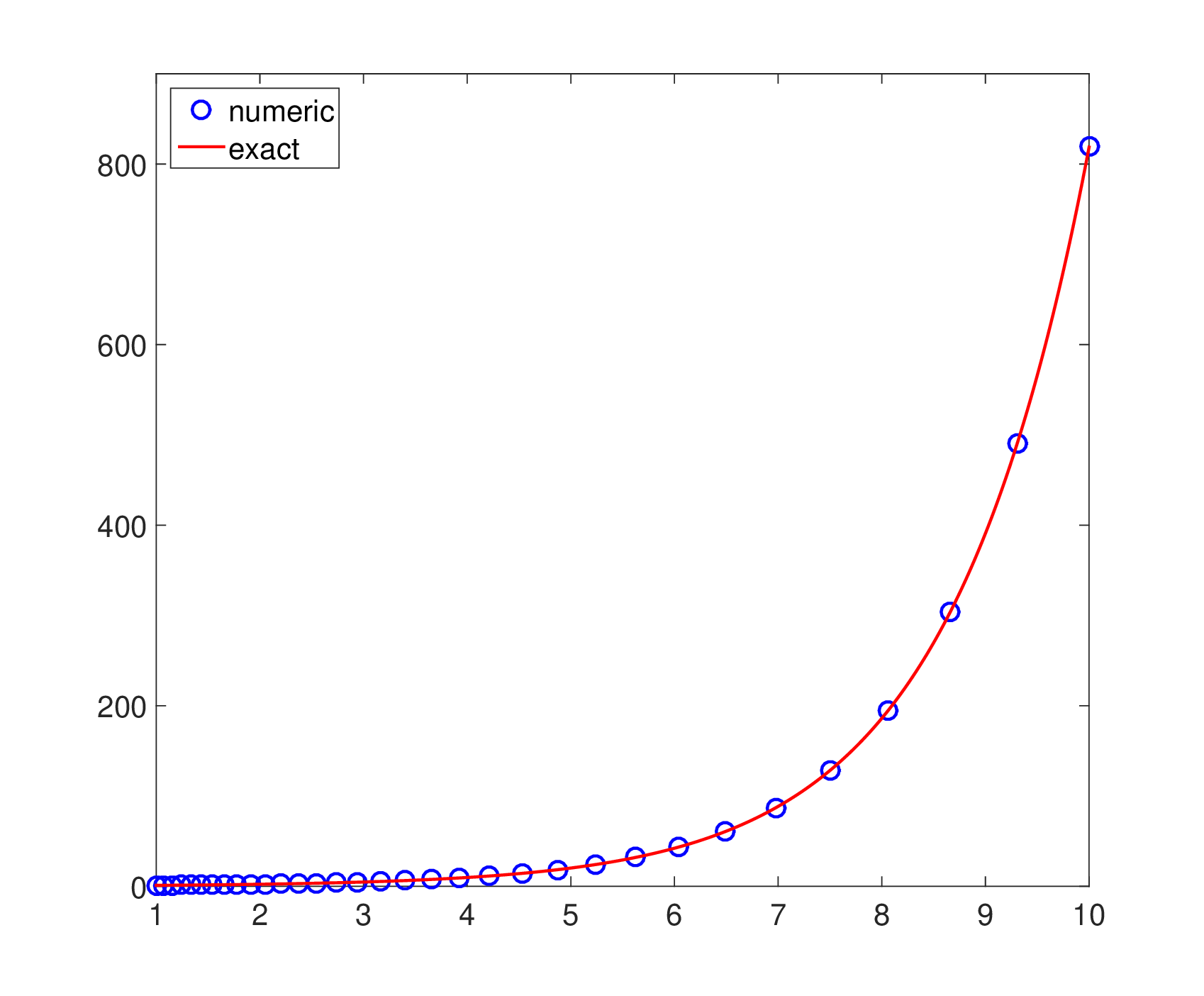}
\caption{Comparison between the numerical (blue circles) and the exact (red solid line) solutions of the system \eqref{Example_2_sys} involving the Caputo-Hadamard derivative operator. The parameters used are $\apa=0.3$, $u_0=1$ and $\lma=1$.}\label{FigHad}
\end{center}
\end{figure}

We use the same finite difference method as in Example \ref{exmp1} to approximate the solution of system \eqref{Example_2_syseqv}, and we obtain an approximation of $u$ by a scaling technique. Figure \ref{FigHad} shows a comparison between the numerical and the exact solutions, and Table \ref{TabHad} lists the $\ell^\infty$ errors and the numerical convergence orders of the scheme for various values of $\apa$. As expected, the theoretical rate of convergence is reached for all the choices of $\apa$.
\begin{table}
\begin{center}
\setlength\extrarowheight{4pt}
\begin{tabular}{@{\extracolsep{\fill}}|c|c|c||c|c||c|c|}
\cline{2-7}
\multicolumn{1}{c|}{}  & \multicolumn{2}{c||}{$\apa=0.7$} & \multicolumn{2}{c||}{$\apa=0.5$} & \multicolumn{2}{c|}{$\apa=0.3$} \\ 
\hline 
$N$ & $\|e\|_{\ell^\infty}$ & Order & $\|e\|_{\ell^\infty}$ & Order & $\|e\|_{\ell^\infty}$ & Order \\
\cline{1-7} 
$2^7$    & $1.0080\ (-02)$ & $--$ & $1.0225\ (-02)$ & $--$ & $4.2867\ (-02)$ & $--$ \\ 
\cline{1-7} 
$2^8$    & $4.1082\ (-03)$ & $1.2949$ & $3.6865\ (-03)$ & $1.4717$ & $1.3865\ (-02)$ & $1.6284$ \\ 
\cline{1-7} 
$2^9$    & $1.6743\ (-03)$ & $1.2950$ & $1.3238\ (-03)$ & $1.4775$ & $4.4766\ (-03)$ & $1.6309$ \\ 
\cline{1-7} 
$2^{10}$ & $6.8176\ (-04)$ & $1.2962$ & $4.7346\ (-04)$ & $1.4834$ & $1.4344\ (-03)$ & $1.6420$ \\
\cline{1-7} 
$2^{11}$    & $2.7739\ (-04)$ & $1.2974$ & $1.6878\ (-04)$ & $1.4881$ & $4.5603\ (-04)$ & $1.6532$ \\ 
\cline{1-7} 
$2^{12}$    & $1.1279\ (-04)$ & $1.2983$ & $6.0023\ (-05)$ & $1.4916$ & $1.4404\ (-04)$ & $1.6627$ \\ 
\hline
\end{tabular}
\vspace*{2mm}
\caption{Errors and convergence orders relative to system \eqref{Example_2_sys} for various values of $\apa$. Notice that the numerical orders of convergence are approximately equal to the theoretical value $2-\apa$.}\label{TabHad}
\end{center}
\end{table}
\begin{rem}
Example 2 illustrates how one can straightforwardly derive a numerical scheme suited to the Hadamard derivative operator from a scheme approximating the standard Caputo derivative while keeping the optimal order of convergence of this latter.
\end{rem}
\subsection{Example 3}
Let $\apa\in (1,2)$ and consider the system 
\begin{align}\label{Example_3_sys}
\left\{
\begin{array}{l}
\Erd{\apa}{0,\sigma,\eta}u(t) = \lma\,t^{\sigma\apa}u(t),  \ \ t \in (0,T]\\ \\
\Erd{\apa-1}{0,\sigma,0}\left(t^{\sigma(1+\eta)}u(t)\right)(0) = a_1, \ \ \lim_{t\to 0}\Iint{2-\apa}{0,\sigma,\apa-2}\left(t^{\sigma(2+\eta)}u(t)\right) = a_2.
\end{array}
\right.
\end{align}
with exact solution
\begin{equation}\label{exact_sol_EK}
u(t) = t^{-\sigma(\apa+\eta)}\left(a_1\,e_\apa^{\lma t^\sigma} + \dfrac{a_2}{\Gamma(\apa-1)t^{\sigma(2-\apa)}}\right),
\end{equation}
where $e_\apa$ is the $\apa$-exponential function given by \eqref{e_alpha}. One may show using Theorem \ref{Theo_EK_equiv} that \eqref{Example_3_sys} is equivalent to the system
\begin{equation*}
\left\{
\begin{array}{l}
\der{\apa}{0}\omega(t) = \lma\,\omega(t), \ \ t\in (0,T^{\sigma}] \\
\der{\apa-1}{0}\omega(0) = a_1, \ \ \lim_{t\to 0}\Iint{2-\apa}{0}\omega(t) = a_2,
\end{array}
\right.
\end{equation*}
where $\omega$ and $u$ are related by $u(t)=t^{-\sigma(\apa+\eta)}\omega(t^\sigma)$. If we set $$z(t) = \omega(t) - \dfrac{a_1}{\Gamma(\apa)}\,t^{\apa-1}- \dfrac{a_2}{\Gamma(\apa-1)}\,t^{\apa-2},$$ then $z$ satisfies
\begin{align}\label{Example_3_syseqv}
\left\{
\begin{array}{l}
\der{\apa}{0}z(t) = \lma\, z(t) + \dfrac{\lma\, a_1}{\Gamma(\apa)}\,t^{\apa-1}, \ \ t\in (0,T^{\sigma}] \\
\der{\apa-1}{0}z(0) = 0, \ \ \lim_{t\to 0}\Iint{2-\apa}{0}z(t) = 0.
\end{array}
\right.
\end{align}
We follow the L2 method introduced in \cite{Old74} in order to obtain an approximation $z_n$ of $z(\bar{t}_n)$, with $\bar{t}_n := \frac{n-1}{N-1}T^{\sigma}$, $1\leq n\leq N$. The approximation of $u(\bar{t}_n^{\ 1/\sigma})$ is then given by 
\begin{equation}\label{num_sol_EK}
u_n = \dfrac{1}{\bar{t}_n^{\ \apa+\eta}}\left(z_n + \dfrac{a_1}{\Gamma(\apa)}\bar{t}_n^{\,\apa-1} + \dfrac{a_2}{\Gamma(\apa-1)}\bar{t}_n^{\,\apa-2}\right), \ \ 2\leq n \leq N.
\end{equation}
\begin{table}
\begin{center}
\setlength\extrarowheight{4pt}
\begin{tabular}{@{\extracolsep{\fill}}|c|c|c||c|c||c|c|}
\cline{2-7}
\multicolumn{1}{c|}{}  & \multicolumn{2}{c||}{$\sigma = 3$} & \multicolumn{2}{c||}{$\sigma = 0.5$} & \multicolumn{2}{c|}{$\sigma = 0.1$} \\ 
\cline{1-7} 
$N$ & $\|e\|_{\ell^\infty}$ & Order & $\|e\|_{\ell^\infty}$ & Order & $\|e\|_{\ell^\infty}$ & Order \\ 
\cline{1-7} 
$2^7$ & $3.7335\ (-02)$ & -- & $3.1528\ (-03)$ & -- & $1.9474\ (-03)$ & -- \\ 
\cline{1-7} 
$2^8$ & $2.1331\ (-02)$ & $0.8076$ & $1.5320\ (-03)$ & $1.0412$ & $9.5087\ (-04)$ & $1.0343$ \\ 
\cline{1-7} 
$2^9$ & $1.2039\ (-02)$ & $0.8253$ & $7.4046\ (-04)$ & $1.0489$ & $4.6223\ (-04)$ & $1.0406$ \\ 
\cline{1-7} 
$2^{10}$ & $6.6692\ (-03)$ & $0.8521$ & $3.5755\ (-04)$ & $1.0503$ & $2.2456\ (-04)$ & $1.0415$ \\ 
\cline{1-7} 
$2^{11}$ & $3.6275\ (-03)$ & $0.8785$ & $1.7291\ (-04)$ & $1.0481$ & $1.0925\ (-04)$ & $1.0395$ \\ 
\cline{1-7} 
$2^{12}$ & $1.9421\ (-03)$ & $0.9014$ & $8.3851\ (-05)$ & $1.0441$ & $5.3274\ (-05)$ & $1.0361$ \\ 
\hline
\end{tabular}
\vspace*{2mm}
\caption{Errors and convergence orders of system \eqref{Example_3_sys} for various values of $\sigma$. The parameters used are: $\apa=1.75$, $a_1=a_2=\lma=1$.  The parameter $\eta$ can be chosen arbitrarily.}\label{TabEK}
\end{center}
\end{table}

\begin{figure}
\begin{center}
\hspace*{-1.7cm}
\begin{tabular}{c}
\includegraphics[scale=0.18]{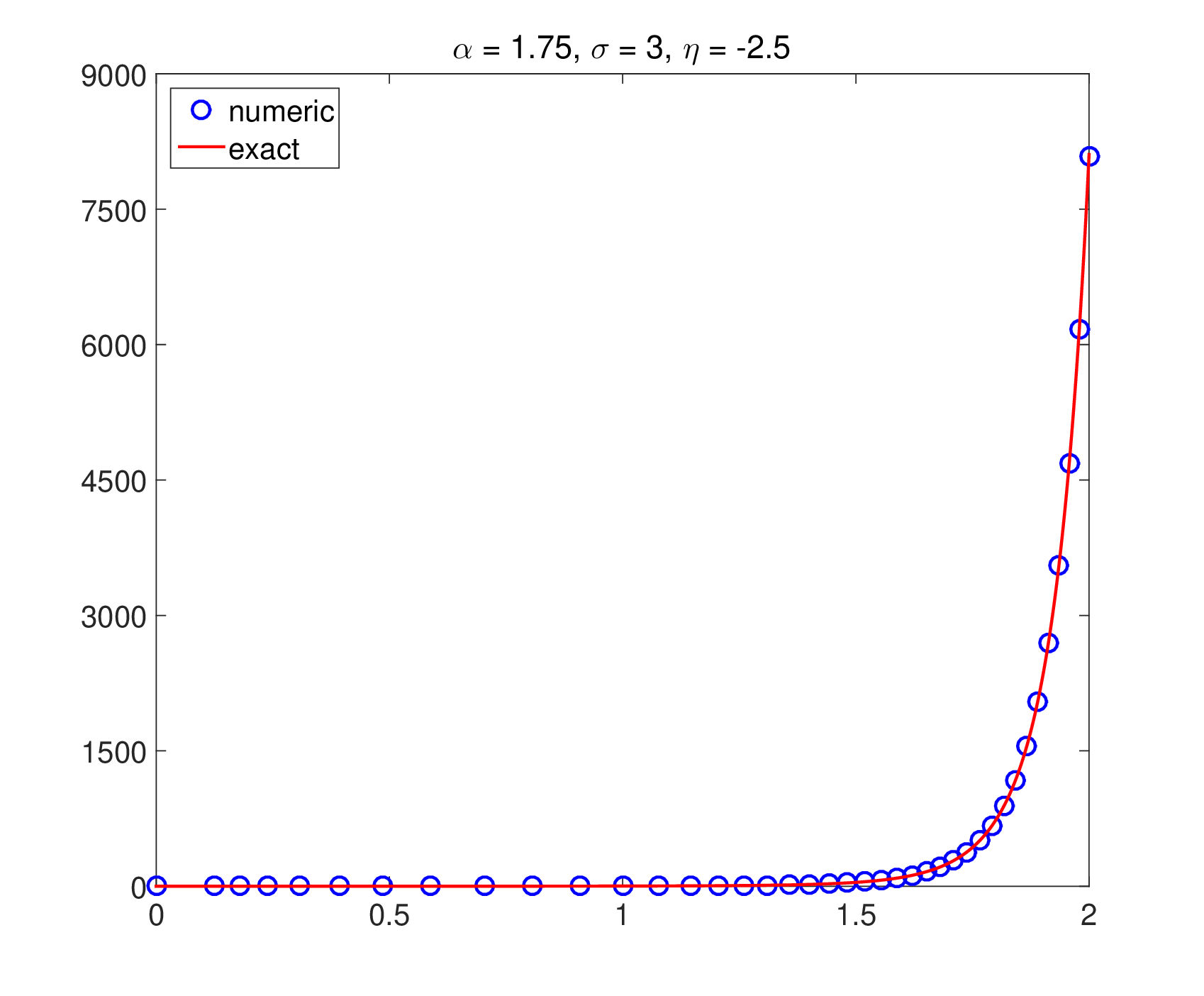} \hspace*{-5mm}\includegraphics[scale=0.18]{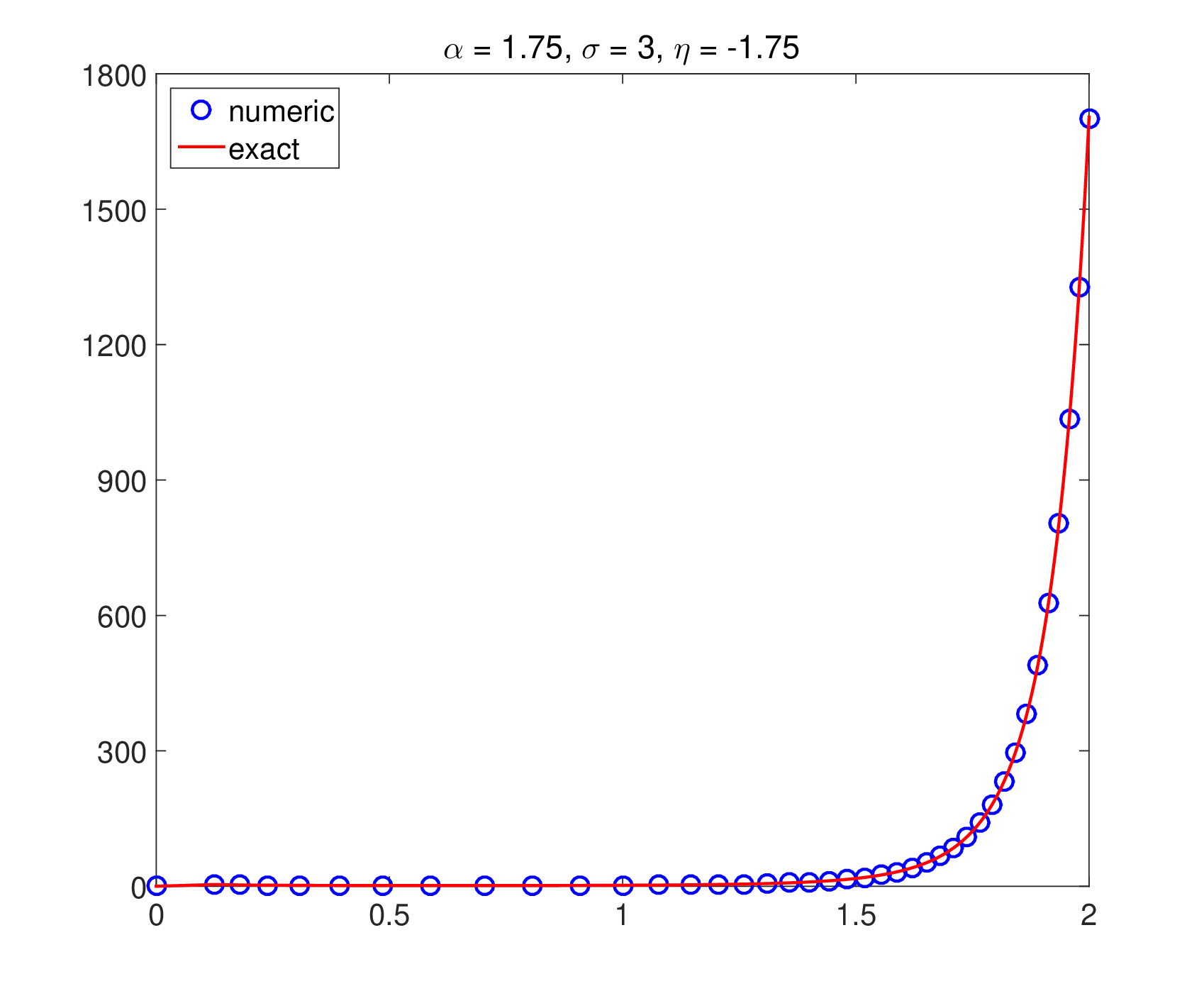} \hspace*{-5mm}\includegraphics[scale=0.18]{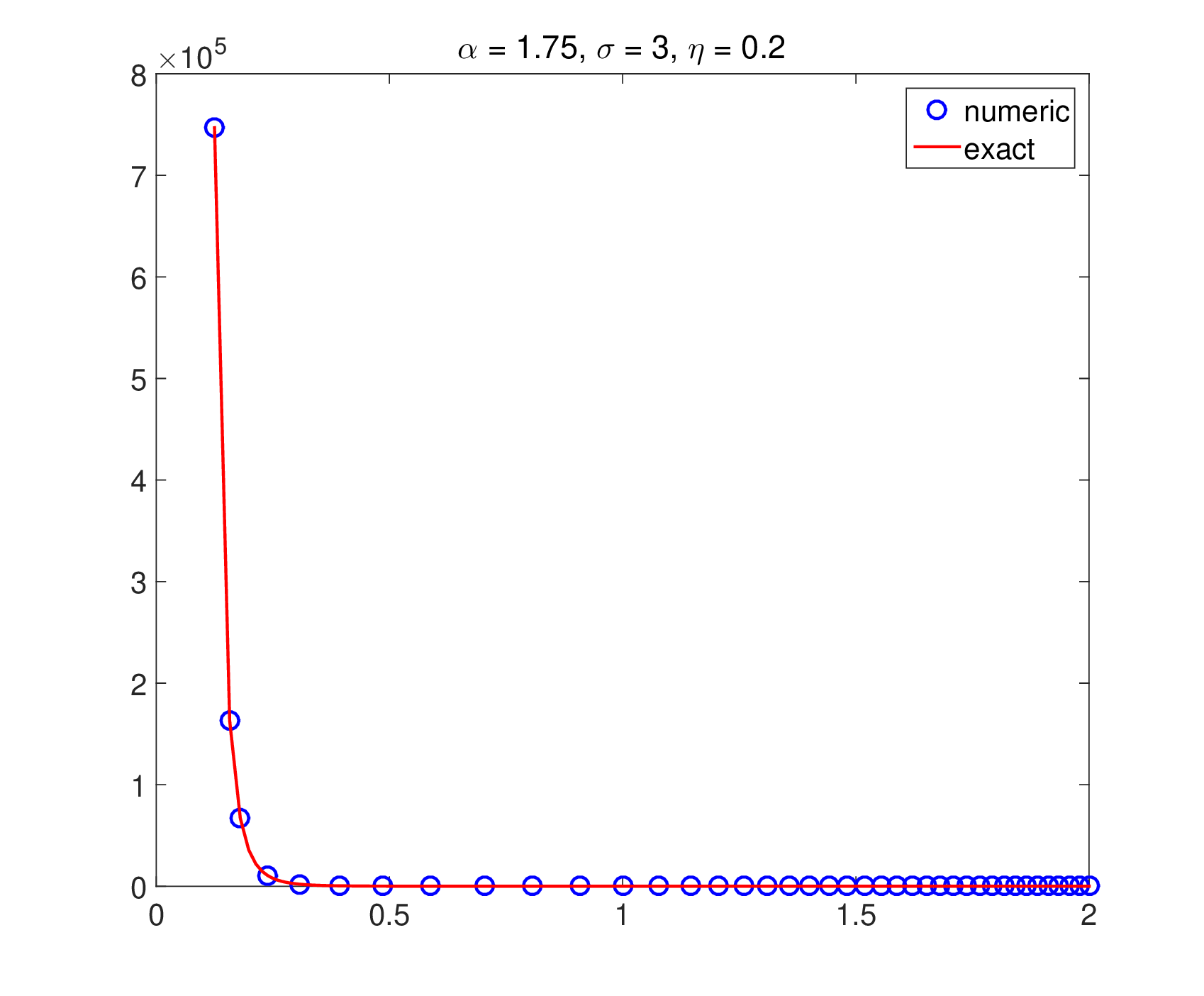}\\
\includegraphics[scale=0.18]{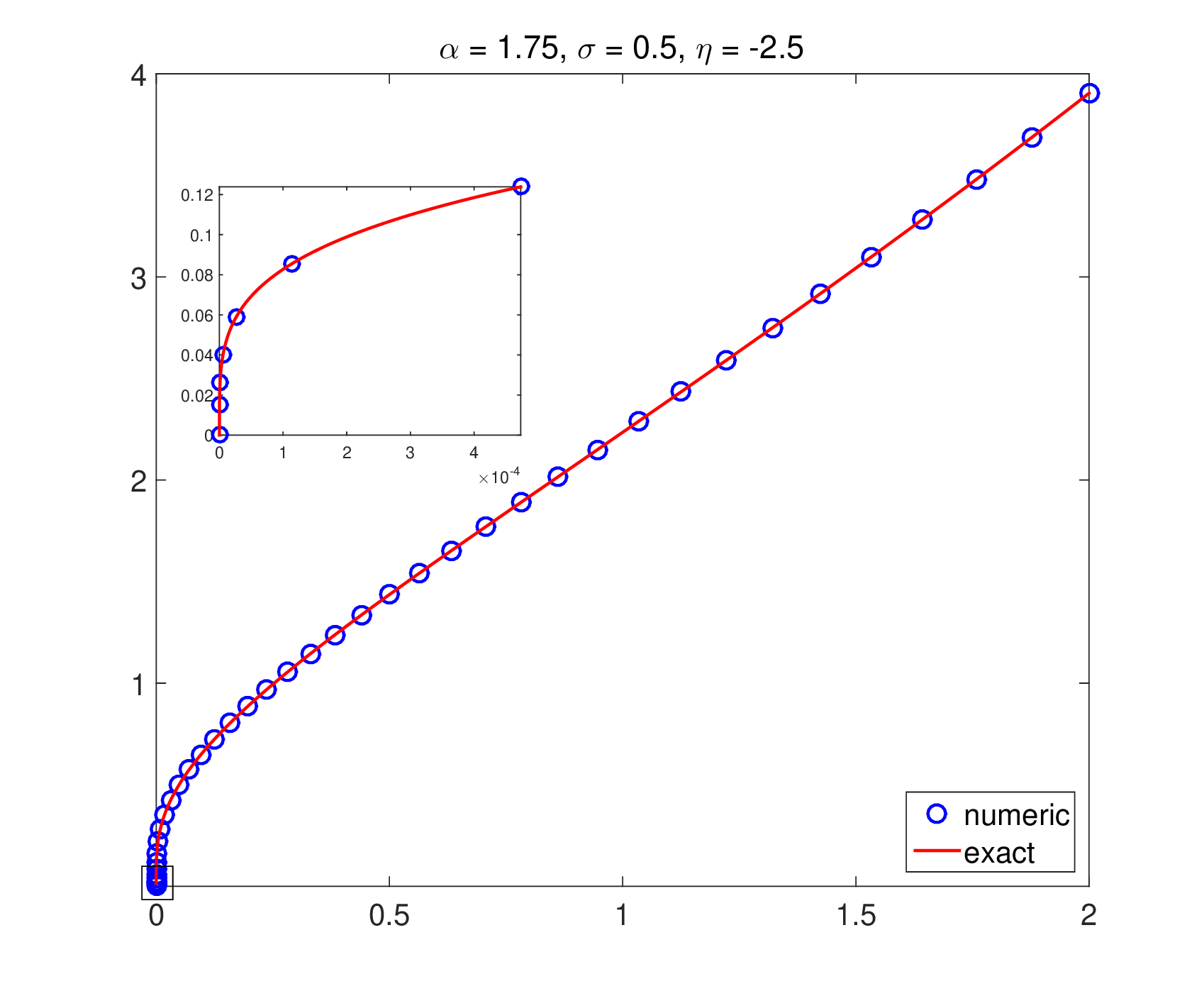} \hspace*{-5mm}\includegraphics[scale=0.18]{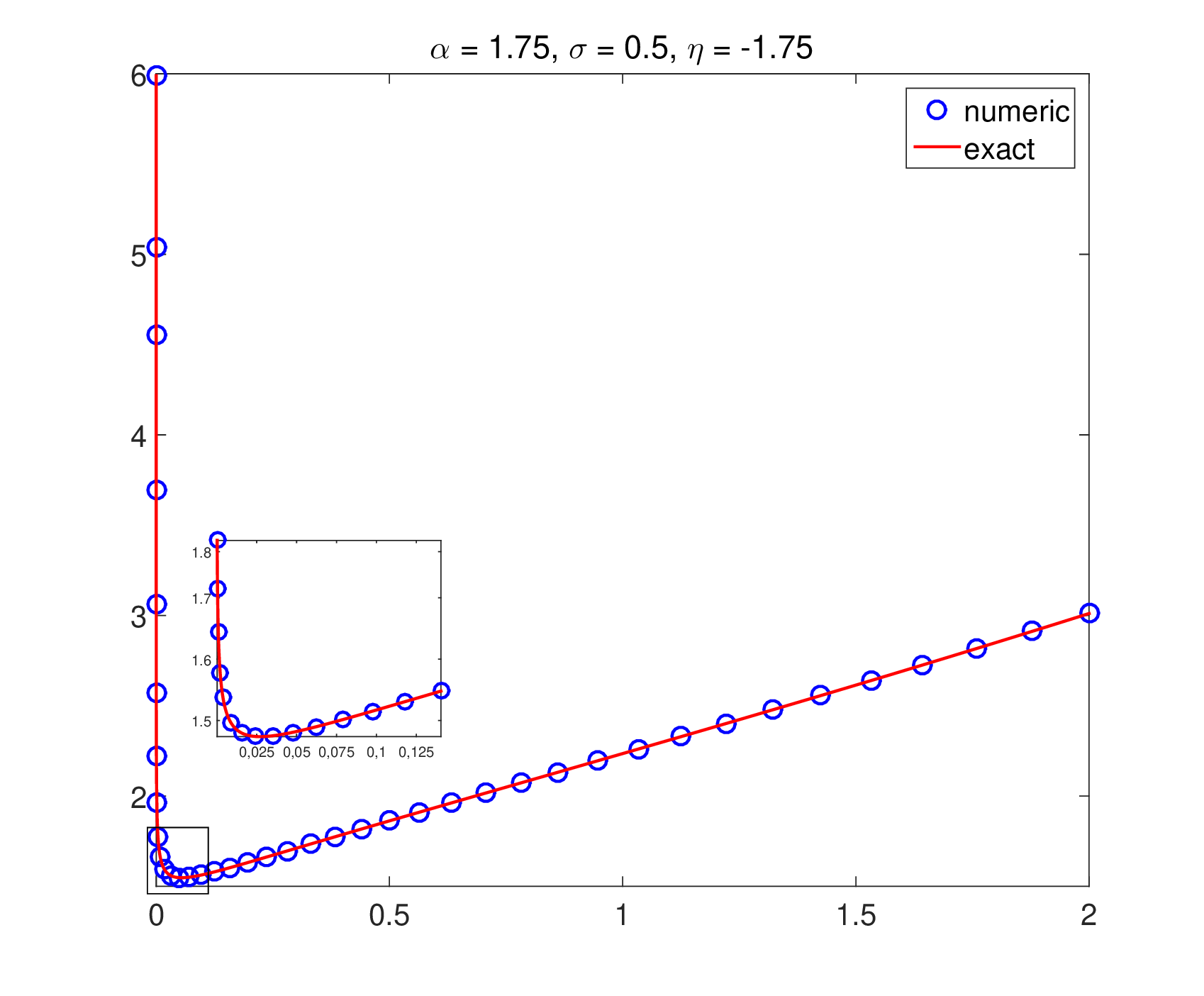} \hspace*{-5mm}\includegraphics[scale=0.18]{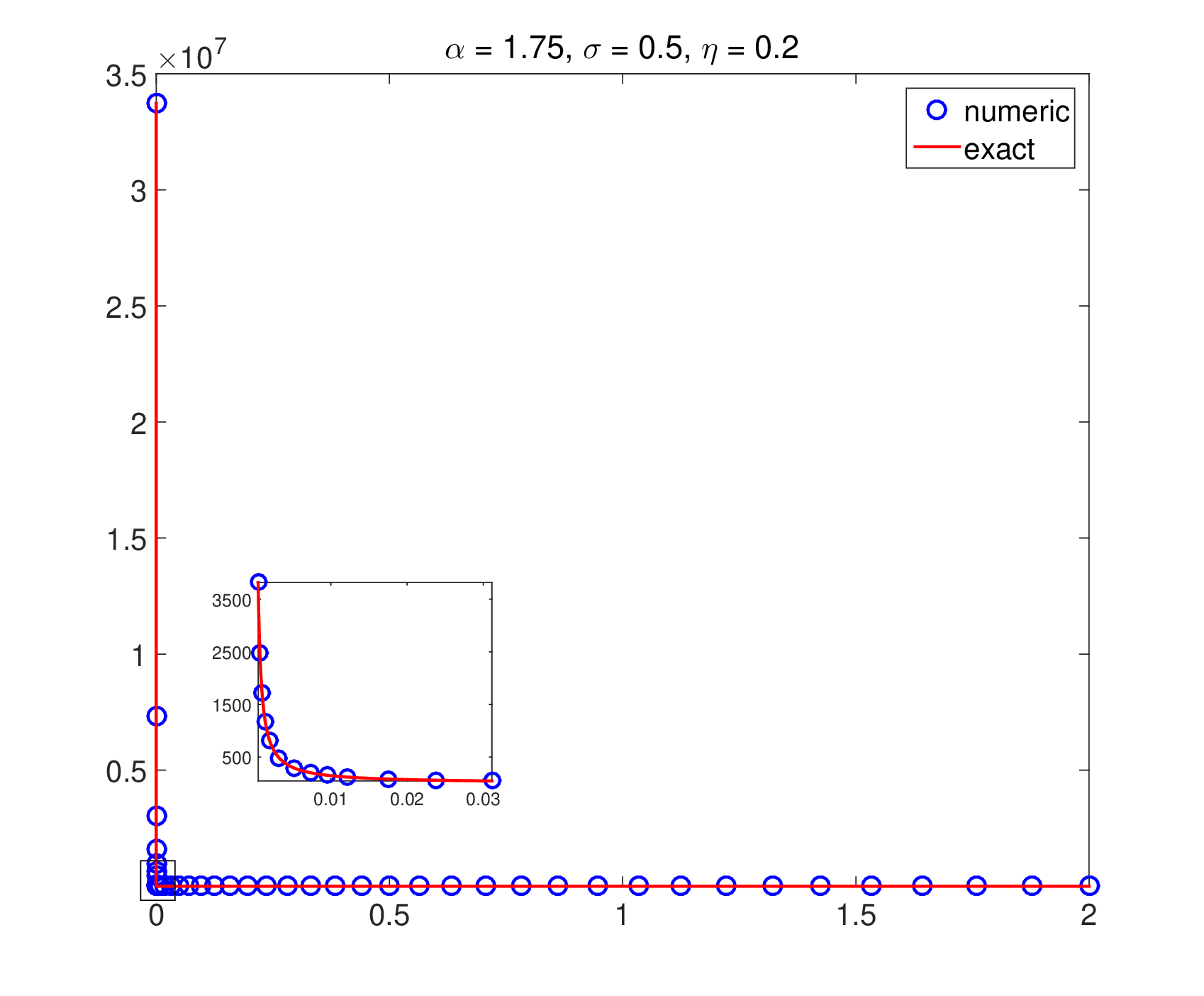}\\
\includegraphics[scale=0.18]{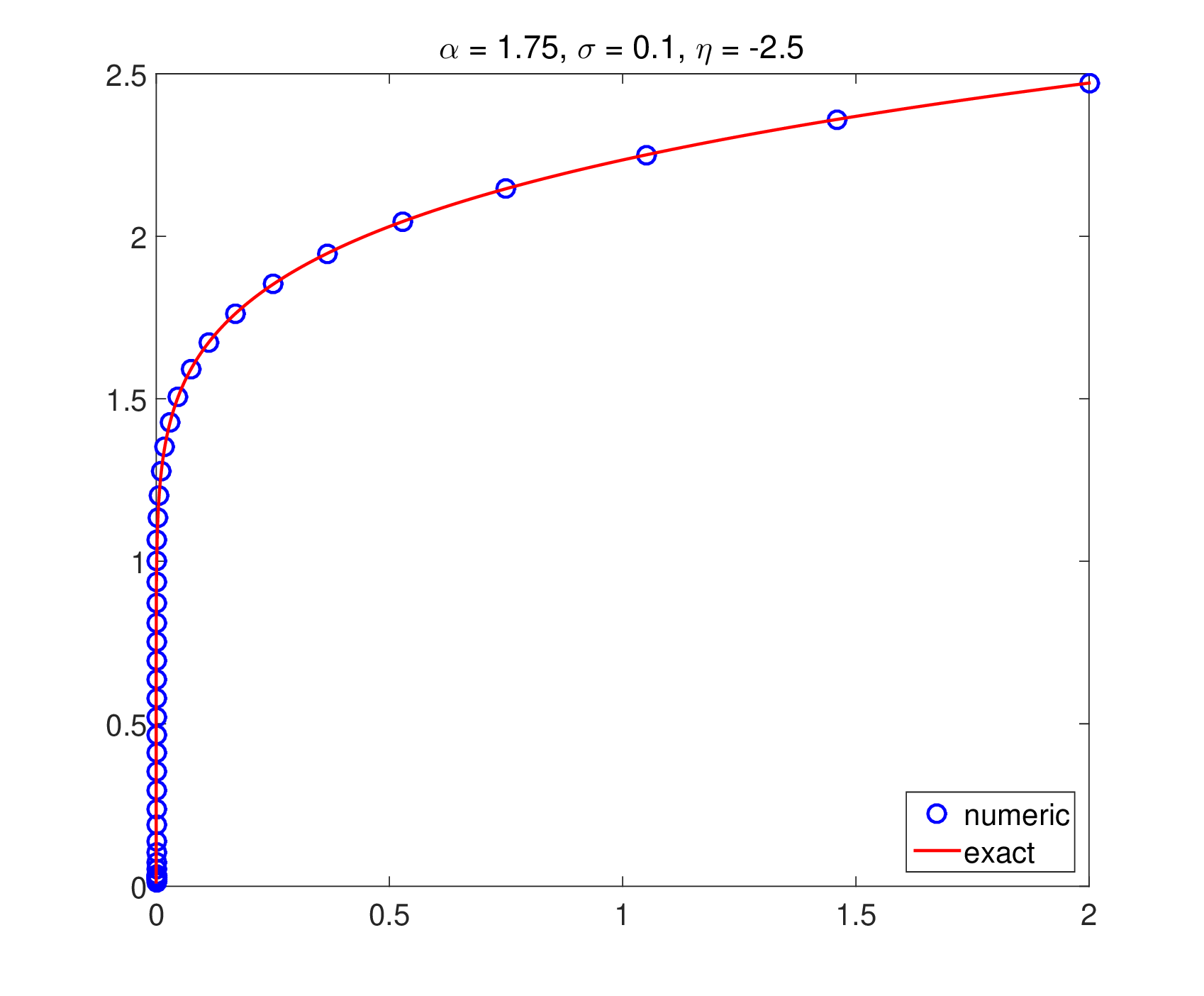} \hspace*{-5mm}\includegraphics[scale=0.18]{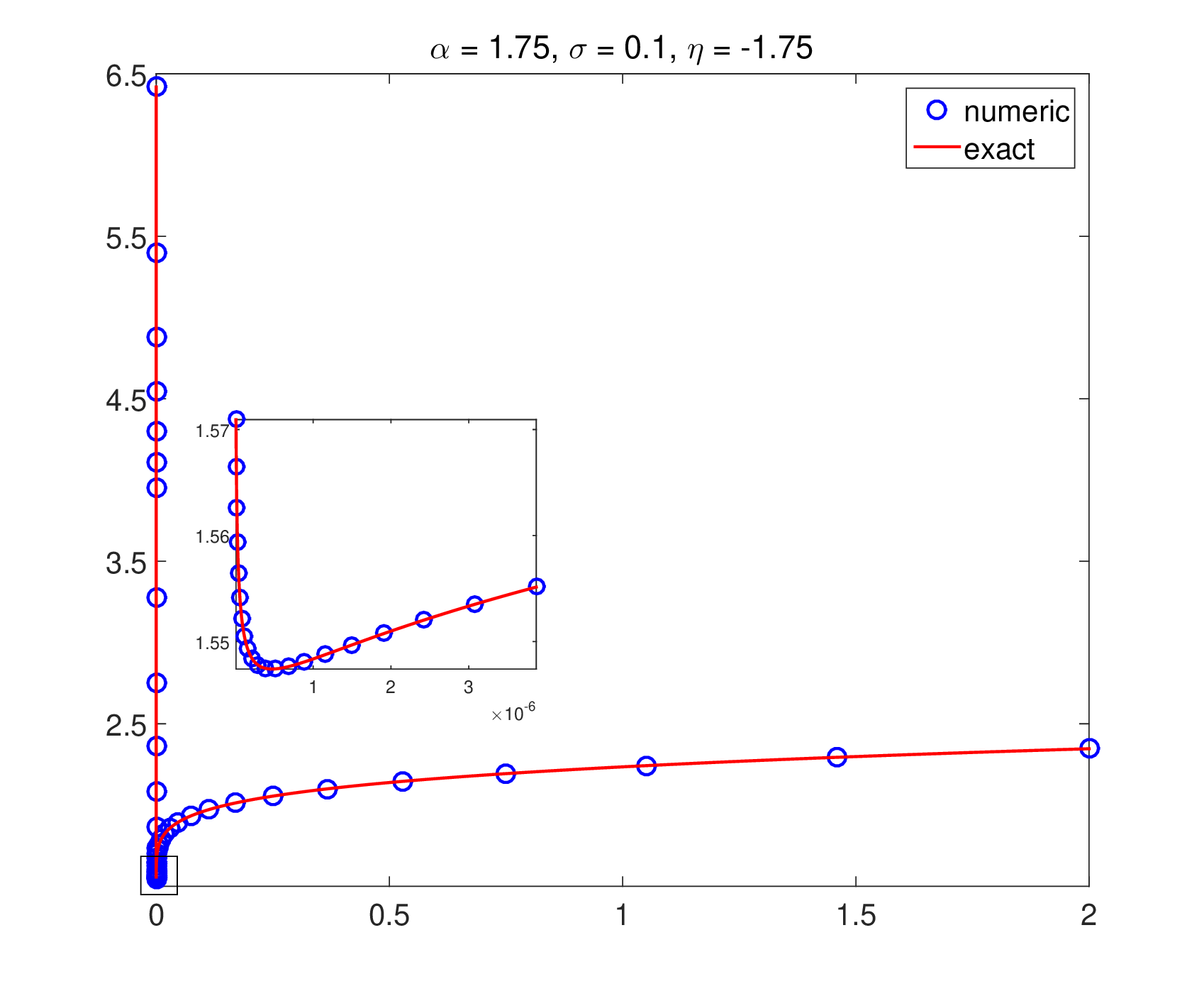} \hspace*{-5mm}\includegraphics[scale=0.18]{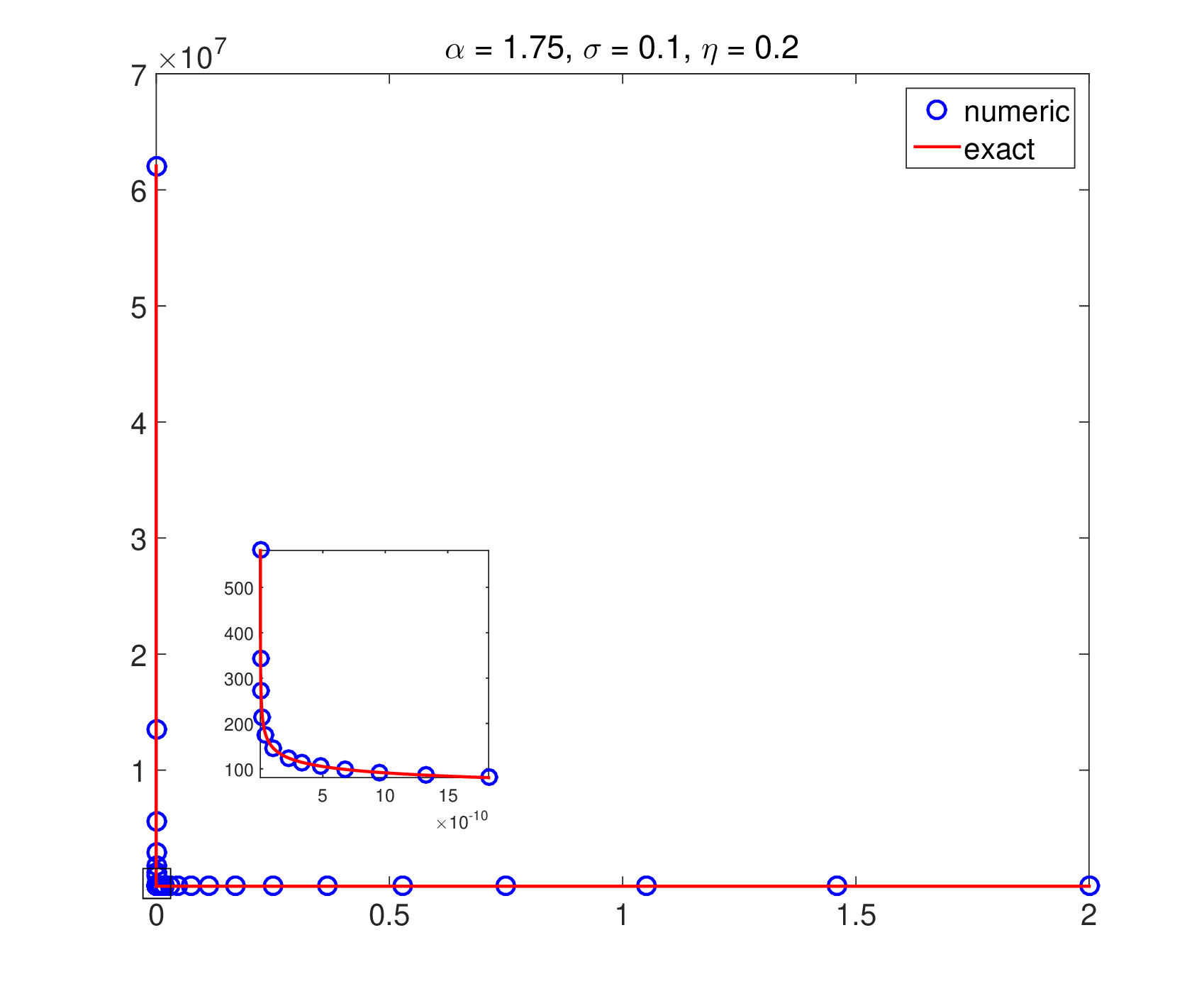}
\end{tabular}
\caption{Numerical solutions (blue circles) versus exact solutions (red solid lines) of the system \eqref{Example_3_sys} involving the Erd\'elyi-Kober derivative operator with $\apa = 1.75$ and for various values of $\sigma$ and $\eta$.}\label{fig_EK}
\end{center}
\end{figure}

Figure \ref{fig_EK} shows a comparison between the numerical solutions \eqref{num_sol_EK} versus the exact solutions \eqref{exact_sol_EK} for various values of $\sigma$ and $\eta$. Though the solutions or their first derivatives might be singular at the origin, one can notice that the numerical solutions fit very well with the exact solutions for all the chosen parameters. In Table \ref{TabEK}, we listed the $\ell^\infty$ errors and the convergence orders for various values of $\sigma$. Obviously, the optimal first order is reached for all the values of $\sigma$ and $\eta$, whether the solutions are regular or not. Moreover, we remarked that only the parameter $\sigma$ is significant to this study while changing the parameter $\eta$ does not impact the numerical rate of convergence. Actually, 
these assertions are expected since the computations are performed using system \eqref{Example_3_syseqv} rather than system \eqref{Example_3_sys}, thus the parameter $\eta$ is not relevant and the meshing nodes only vary when the value of $\sigma$ changes. This confirms the robustness of our approach which allows to accurately approach the solutions of systems involving the Erd\'elyi-Kober operators in a simple and general framework.

\section*{Conclusion and perspectives}
This paper is concerned with the study of the fractional integrals and derivatives with respect to another function. By establishing a one-to-one correspondence between the fractional operators with respect to another function and the standard Riemann-Liouville or Caputo fractional operators in scaled axes (see Corollary \ref{CoroEq_RL}, Theorem \ref{TheoEqCap}, and  Theorem \ref{Theo_EK_equiv}), we proved several results related to the fractional calculus in appropriate functional spaces. We also showed that any numerical scheme for the RL or Caputo operators can adequately be used to approach the solutions of systems involving the fractional operators with respect to another function. Though the solutions might be singular, the approximated solutions are generated on graded meshes and the convergence orders of the numerical schemes remain optimal (i.e. do not depend on the scaling function). Our approach can be applied to any fractional operator that can be expressed in terms of fractional derivative with respect to another function, such as the Hadamard derivative \cite{Kilbas06}, the generalized derivative \cite{Katugampola14,Zen17}, the Erd\'elyi-Kober derivative \cite{Kiryakova94,Kiryakova13}, and so on, depending on the choice of the function $\psi$, without any additional computational costs. The efficiency of the proposed method is highlighted throughout several examples and numerical tests. As perspectives of this work, one could extend the proposed method to study other fractional operators such as the $\psi$-Hilfer operators \cite{Hilfer00} or the fractional operators with respect to two functions \cite{Agr12}, but also the fractional boundary value problems \cite{Abd19}. This will be the subject of a future work.

\begin{appendices}
\section{}\label{SectionAppA}
\begin{lemma}\label{lem_SectionAppA}
Let $n\in \N$ and let $\psi\in C^n[a,b]$ be a monotonous function such that $\psi'(t)\neq 0$ for all $t\in [a,b]$. Then we have the following embedding
\begin{equation*}
C^n[a,b] \subset AC^n_\psi[a,b] \subset C^{n-1}[a,b] \subset \dots \subset C^1[a,b] \subset AC_\psi[a,b] \subset C[a,b]
\end{equation*}
where $C^n[a,b]$ denotes the set of continuously differentiable functions up to order $n$.
\end{lemma}
\begin{proof}
Let $u\in C^n[a,b]$. Since $\psi' \in C^{n-1}[a,b]$ with $\psi'(t)\neq 0$ for all $t\in [a,b]$, then the function $\varrho u = \frac{u'}{\psi'}\in C^{n-1}[a,b]$. Applying $\varrho$ again yields $\varrho^2 u = \varrho(\varrho u)\in C^{n-2}[a,b]$. Applying $\varrho$ recursively yields $\varrho^k u \in C^{n-k}[a,b]$ for all $0\leq k \leq n-1$. In particular $\varrho^{n-1} u \in C^1[a,b] \subset AC[a,b]$, and hence $u\in AC^n_\psi[a,b]$.\\

Now we prove that $AC^n_\psi[a,b] \subset C^{n-1}[a,b]$. Let $u\in AC^n_\psi[a,b]$, hence $\varrho^{n-1} u \in AC[a,b] \subset C[a,b]$, i.e. $\varrho(\varrho^{n-2} u) \in C[a,b]$. Since $\psi'\in C[a,b]$ then $(\varrho^{n-2} u)' \in C[a,b]$, or equivalently $\varrho^{n-2} u \in C^1[a,b]$. Proceeding similarly, one obtain $\varrho^{n-k} u \in C^{k-1}[a,b]$ for all $1 \leq k \leq n$, and thus $u \in C^{n-1}[a,b]$, which ends the proof.
\end{proof}
\begin{rem}
Lemma \ref{lem_SectionAppA} is no more valid if one suppose $\psi'(t)\neq 0$ for all $t\in (a,b]$ (instead of $\psi'(t)\neq 0$ for all $t\in [a,b]$). For instance, one may consider $n=2$, $\psi(t)=\frac{t^\rho}{\rho}$ and $u=id$ in $[0,1]$. Then $\forall\ t\in (0,1)$
\begin{align*}
\left(\varrho u\right)'(t)=\left(\dfrac{du}{d\psi}\right)'(t)=\Big{(}t^{1-\rho}\,u'(t)\Big{)}' = (1-\rho)t^{-\rho}\not\in L^1(0,1) \ \text{ if }\ \rho> 1,
\end{align*}
and hence $C^2[0,1]\not\subset AC^2_\psi[0,1]$ for $\rho > 1$.
\end{rem}
\section{}\label{SectionAppB}
\begin{lemma}
Let $k\in \N$ and suppose $f$ is $k$ times derivable on $\R_+^*$. Then
\begin{align}\label{gamma_kf}
\left(\varrho^ku\right)(t) = \dfrac{1}{\sigma^k}\sum_{j=1}^k \lma_{j,k}\,t^{j-k\sigma}u^{(j)}(t)
\end{align}
where $(\lma_{j,k})_{\substack{0\leq j \leq m \\ 1\leq k \leq m}}$ is the sequence defined by
\begin{align*}
\lma_{j,k} = \left\{
\begin{array}{ll}
0 & \ \text{ if }\ \ j=0 \ \text{ or }\ j > k \\
1 & \ \text{ if }\ \ j = k \\
\lma_{j-1,k-1} + \big{(}j-(k-1)\sigma\big{)}\,\lma_{j,k-1} & \ \text{ if }\ \ 1\leq j < k.
\end{array}
\right.
\end{align*}
\end{lemma}
\begin{proof}
We prove \eqref{gamma_kf} by induction. The result is trivial for $k=1$. Assume \eqref{gamma_kf} holds true, then
\begin{align*}
\left(\varrho^{k+1}u\right)(t) & = \varrho\left(\varrho^ku\right)(t) \\
& = \dfrac{1}{\sigma}t^{1-\sigma}\,\dfrac{d}{dt}\left(\dfrac{1}{\sigma^k}\sum_{j=1}^k \lma_{j,k}\,t^{j-k\sigma}u^{(j)}(t)\right)\\
& = \dfrac{1}{\sigma^{k+1}}t^{1-\sigma}\left(\sum_{j=1}^k (j-k\sigma)\,\lma_{j,k}\,t^{j-k\sigma-1}u^{(j)}(t) + \sum_{j=1}^k \lma_{j,k}\,t^{j-k\sigma}u^{(j+1)}(t)\right)\\
& = \dfrac{1}{\sigma^{k+1}}\left(\sum_{j=1}^k (j-k\sigma)\,\lma_{j,k}\,t^{j-(k+1)\sigma}u^{(j)}(t) + \sum_{j=2}^{k+1} \lma_{j-1,k}\,t^{j-(k+1)\sigma}u^{(j)}(t)\right)\\
& = \dfrac{1}{\sigma^{k+1}}\left[\sum_{j=1}^{k+1} \left(\lma_{j-1,k} + (j-k\sigma)\,\lma_{j,k}\right)t^{j-(k+1)\sigma}u^{(j)}(t) \ - \underbrace{\lma_{0,k}}_{=\ 0}\,t^{1-(k+1)\sigma}u'(t)\right.\\
& \qquad\qquad\quad \left. -(k+1-k\sigma)\,\underbrace{\lma_{k+1,k}}_{=\ 0}\,t^{(k+1)(1-\sigma)}u^{(k+1)}(t)\right] \\
& = \dfrac{1}{\sigma^{k+1}}\sum_{j=1}^{k+1} \lma_{j,k+1}\,t^{j-(k+1)\sigma}u^{(j)}(t).
\end{align*}
\end{proof}

\end{appendices}



\begin{thebibliography}{1}
%
\bibitem{Abd19}{Abdo, M.S., Panchal, S.K. and Saeed, A.M., Fractional boundary value problem with $\psi$-Caputo fractional derivative, Proc. Math. Sci., \textbf{129}, 65 (2019)}
%
\bibitem{Agr12}{Agrawal, O.P., Some generalized fractional calculus operators and their applications in integral equations, Fract. Calc. Appl. Anal., \textbf{15}(4), 700--711 (2012)}
%
\bibitem{Alm17}{Almeida, R., A Caputo fractional derivative of a function with respect to another function, Commun. Nonlinear. Sci. Numer. Simulat., \textbf{44}, 460--481 (2017)}
%
\bibitem{Alm18}{Almeida, R., Malinowska, A.B. and Monteiro, M.T., Fractional differential equations with a Caputo derivative with respect to a Kernel function and their applications, Math. Meth. Appl. Sci., {\bf 41}(1), 336--352 (2018)}
%
\bibitem{Alm19}{Almeida, R., Further properties of Osler's generalized fractional integrals and derivatives with respect to another function, Rocky Mountain, \textbf{49}(8), 2459--2493, 2019}
%
\bibitem{Bal17}{Baleanu, D., Wu, G.C. and Zeng, S., Chaos analysis and asymptotic stability of generalized Caputo fractional differential equations, Chaos, Solitons \& Fractals \textbf{102}, (2017) 99--105.}
%
\bibitem{Benjemaa18}{Benjemaa, M.: Taylor's formula involving generalized fractional derivatives, Appl. Math. Comput., \textbf{335}, 182--195 (2018)}
%
\bibitem{Cai20}{Cai, M. and Li, C., Numerical Approaches to Fractional Integrals and Derivatives: A Review. Mathematics, {\bf 8}(1), 43, (2020)}
%
\bibitem{Cao2003}{Cao, Y., Herdman, T. and Xu, Y., A hybrid collocation method for Volterra integral equations with weakly singular kernels, Siam J. Numer. Anal., \textbf{41}(1), 364--381 (2003)}
%
\bibitem{Cap71}{Caputo, M. and Mainardi, F., A new dissipation model based on memory mechanism, Pure and Applied Geophysics, \textbf{91}(8), 134--147 (1971)}
%
\bibitem{Cap92}{Caputo, M., Lectures on Seismology and Rheological Tectonics, Univ. degli studi di Roma "La Sapienza", 1992}
%
\bibitem{Col18}{Colombaro, I., Garra, R., Giusti, A. and Mainardi, F., Scott-Blair models with time varying viscosity, Appl. Math. Lett., \textbf{86}, 57--63 (2018)}
%
\bibitem{Diethelm97}{Diethelm, K., An algorithm for the numerical solution of differential equations of fractional order, Elec. Trans. Num. Anal., \textbf{5}, 1--6 (1997)}
%
\bibitem{Diethelm10}{Diethelm, K.: The analysis of fractional differential equations, An Application-Oriented Exposition Using Differential Operators of Caputo Type, Lecture Notes in Mathematics, Springer-Verlag Berlin Heidelberg, (2010)}
%
\bibitem{Erd64}{Erd\'elyi, A., An integral equation involving Legendre functions, SIAM J. Appl. Math., {\bf 12}, 15-30 (1964)}
%
\bibitem{Erd65}{Erd\'elyi, A., Axially symmetric potentials and fractional integration, SIAM J. Appl. Math., {\bf 13}, 216-228 (1965)}
%
\bibitem{Ford2015}{Ford, N.J., Morgado, M.L. and Rebelo, M., A nonpolynomial collocation method for fractional terminal value problems, J. Comput. Appl. Math., \textbf{275}, 392--402 (2015)}
%
\bibitem{Fukunaga15}{Fukunaga, M. and Shimizu, N., Fractional derivative constitutive models for finite deformation of viscoelastic materials, J. Comput. Nonlinear Dyn. \textbf{10}(6), (2015) 061002.}
%
\bibitem{Gar18}{Garra, R., Giusti, A. and Mainardi, F., The fractional Dodson diffusion equation: a new approach, Ricerche Mat. \textbf{67}, 899--909 (2018)}
%
\bibitem{Herrmann14}{Herrmann, R., Fractional Calculus: An Introduction for Physicists. World Scientific Publishing, Singapore, 2\ups{nd} edition, 2014.}

\bibitem{Hilfer00}{Hilfer, R., Applications Of Fractional Calculus In Physics. World Scientific Publishing, Singapore, 2000.}
%
\bibitem{Jar20}{Jarad, F. and Abdeljawad, Th., Generalized fractional derivatives and Laplace transform, Disc. Cont. Dyn. Sys. Series S, \textbf{13}(3), 709--722 (2020)}
%
\bibitem{Katugampola14}{Katugampola, U.N., A New approach to generalized fractional derivatives, Bull. Math. Anal. Appl., \textbf{6}(4), 1--15 (2014)}
%
\bibitem{Kilbas06}{Kilbas, A.A., Srivastava, H.M. and Trujillo, J.J., Theory and Applications of Fractional Differential Equations, North-Holland mathematics studies, {\bf 204}. Amsterdam, Elsevier Science B.V., (2006)}
%
\bibitem{Kiryakova94}{Kiryakova, V., Generalized Fractional Calculus and Applications. Longman
\& J. Wiley, Harlow, New York, 1994.}
%
\bibitem{Kiryakova13}{Kiryakova, V. and Luchko, Y., Riemann-Liouville and Caputo type multiple Erd\'elyi-Kober operators. Central Europ. J. Phys., \textbf{11}(10), 1314--1336 (2013)}
%
\bibitem{Lak09}{Lakshmikantham, V., Leela, S. and Vasundhara, D.J., Theory of fractional dynamic systems. Cambridge Scientific Publishers, 2009.}
%
\bibitem{Lin2007}{Lin, Y. and Xu, C., Finite difference/spectral approximations for the time-fractional diffusion equation, J. Comput. Phys., \textbf{225}, 1533--1552 (2007)}
%
\bibitem{Lubich83}{Lubich, Ch., Runge-Kutta theory for Volterra and Abel integral equations of the second kind, Math. Comput., \textbf{41}(163), 87--102 (1983)}
%
\bibitem{Machado17}{Machado, J.A., Fractional Calculus: Fundamentals and Applications, Acou. Vibr. Mech. Struct., Springer Proceedings in Physics, \textbf{198}, (2017) 3--11.}
%
\bibitem{Miller71}{Miller, R.K. and Feldstein, A., Smoothness of solutions of Volterra integral equations with weakly singular kernels, Siam J. Math. Anal., \textbf{2}(2), 242--258 (1971)}
%
\bibitem{Mil93}{Miller, K.S. and Ross, B., An introduction to the fractional calculus and differential equations. John Wiley, New York, 1993.}
%
\bibitem{Mokh20}{Mokhtari, R. and Mostajeran, F., A High Order Formula to Approximate the Caputo Fractional Derivative, Comm. Appl. Math. Comput., \textbf{2}, 1--29 (2020)}
%
\bibitem{Odi07}{Odibat, Z. and Shawagfeh, N.T., Generalized Taylor's formula, Appl. Math. Comput., \textbf{186}, 286--293 (2007)} 
%
\bibitem{Odi20}{Odibat, Z. and Baleanu, D., Numerical simulation of initial value problems with generalized Caputo-type fractional derivatives, Appl. Num. Math., \textbf{156}, 94--105 (2020)} 
%
\bibitem{Old74}{Oldham, K.B. and Spanier, J., The Fractional Calculus, Academic Press, New York, 1974}
%
\bibitem{Osl70a}{Osler, T.J., The fractional derivative of a composite function, SIAM. J. Math. Anal., {\bf 1}(2), 288--293 (1970)}
%
\bibitem{Osl70b}{Osler, T.J., Leibniz rule for fractional derivatives generalized and an application to infinite series, SIAM. J. Appl. Math., {\bf 18}(3), 658--674 (1970)}
%
\bibitem{Osl71}{Osler, T.J., Taylor's series generalized for fractional derivatives and applications, SIAM. J. Math. Anal., {\bf 2}(1), 37--48 (1971)}
%
\bibitem{Osl72}{Osler, T.J., A further extension of the Leibniz rule to fractional
derivatives and its relation to Parseval's formula, SIAM. J. Math. Anal., {\bf 3}(1), 1--16 (1972)}
%
\bibitem{Osl73}{Osler, T.J., A correction to Leibniz rule for fractional derivatives, SIAM. J. Math. Anal., {\bf 4}(3), 456--459 (1973)}
%
\bibitem{Podlubny99}{Podlubny, I., Fractional Differential Equations. Academic Press, San Diego, California, 1999}
%
\bibitem{Pol08}{Polyanin, A.D. and Manzhirov, A.V., Handbook of Integral Equations.
2\ups{nd} Ed., Chapman \& Hall/CRC, Boca Raton, FL 2008}
%
\bibitem{Rice1969}{Rice, J., On the degree of convergence of nonlinear spline approximations, in Approximation with Special Emphasis on Spline Functions, I.J. Schoenberg, ed., Academic Press, NY, 349--365 (1969)}
%
\bibitem{Samko93}{Samko, S.G., Kilbas A.A. and Marichev, O.I., Fractional Integrals and Derivatives: Theory and Applications, Gordon and Breach Science Publishers, Switzerland, 1993}
%
\bibitem{Sne75}{Sneddon, I.N., The use in mathematical physics of Erd\'elyi-Kober operators and of some of their generalizations. In: Ross B. (eds) Fractional Calculus and Its Applications. Lecture Notes in Mathematics, \textbf{457}. Springer, Berlin, Heidelberg. 1975.}
%
\bibitem{Sne79}{Sneddon, I.N., The Use of Operators of Fractional Integration in Applied Mathematics, PWN - Polish Sci. Publishers, Warszawa-Poznan, 1979}
%
\bibitem{Sou18}{Sousa, J.C. and Oliveira, E.C., On the $\psi$-Hilfer fractional derivative, Comm. Nonlin. Sci. Numer. Simul., {\bf 60}, 72--91 (2018)}
%
\bibitem{Yan20}{Yang, Y. and Ji, D., Properties of positive solutions for a fractional boundary value problem involving fractional derivative with respect to another function, AIMS Mathematics, {\bf 5}(6), 7359--7371 (2020)}
%
\bibitem{Yuf13}{Yufeng, X., Zhimin, H., Qinwu, X., Numerical solutions of fractional advection-diffusion equations with a kind of new generalized fractional derivative, International Journal of Computer Mathematics, \textbf{91}, 37--41 (2013)}
%
\bibitem{Yufe13}{Yufeng, X., Agrawal, O.P, Numerical solutions and analysis of diffusion for new generalized fractional Burgers equation, Fractional calculus and applied analysis, \textbf{16}, 709--739 (2013) }
%
\bibitem{Yufen13}{Yufeng, X., Zhimin, H., Qinwu, X., Numerical and analytical solutions of new generalized fractional diffusion equation, Computers and Mathematics with applications \textbf{66}, 2019--2029 (2013)}
%
\bibitem{Zen17}{Zeng, S., Baleanu, D., Bai, Y. and Wu, G., Fractional differential equations of Caputo-Katugampola type and numerical solutions, Appl. Math. Comput., \textbf{315}, 549--554 (2017)}
\end{thebibliography}
\end{document}